\documentclass[reqno,10pt,centertags]{amsart}
\usepackage{amsmath,amsthm,amscd,amssymb,latexsym,upref}

\newcommand{\bbN}{{\mathbb{N}}}
\newcommand{\bbR}{{\mathbb{R}}}

\newcommand{\bbZ}{{\mathbb{Z}}}
\newcommand{\bbC}{{\mathbb{C}}}

\newcommand{\calB}{{\mathcal B}}

\newcommand{\calD}{{\mathcal D}}

\newcommand{\calH}{{\mathcal H}}

\newcommand{\calK}{{\mathcal K}}

\newcommand{\calM}{{\mathcal M}}

\newcommand{\dott}{\,\cdot\,}
\newcommand{\no}{\nonumber}
\newcommand{\lb}{\label}
\newcommand{\f}{\frac}
\newcommand{\ul}{\underline}
\newcommand{\ol}{\overline}
\newcommand{\ti}{\tilde  }
\newcommand{\wti}{\widetilde  }

\newcommand{\loc}{\text{\rm{loc}}}
\newcommand{\spec}{\text{\rm{spec}}}
\newcommand{\rank}{\text{\rm{rank}}}

\newcommand{\dom}{\text{\rm{dom}}}

\newcommand{\AC}{\text{\rm{AC}}}
\newcommand{\bi}{\bibitem}

\newcommand{\beq}{\begin{equation}}
\newcommand{\eeq}{\end{equation}}
\newcommand{\ba}{\begin{align}}
\newcommand{\ea}{\end{align}}

\renewcommand{\Re}{\text{\rm Re}}
\renewcommand{\Im}{\text{\rm Im}}

\allowdisplaybreaks
\numberwithin{equation}{section}

\newtheorem{theorem}{Theorem}[section]
\newtheorem{lemma}[theorem]{Lemma}
\newtheorem{corollary}[theorem]{Corollary}
\newtheorem{hypothesis}[theorem]{Hypothesis}
\theoremstyle{definition}
\newtheorem{definition}[theorem]{Definition}

\theoremstyle{remark}
\newtheorem{remark}[theorem]{Remark}

\begin{document}
\title[Uniqueness for Schr\"odinger, 
Jacobi, and Dirac-Type 
Operators]{Uniqueness Results for Matrix-Valued Schr\"odinger, 
Jacobi, and Dirac-Type Operators}
\author[Gesztesy, Kiselev, and
Makarov]{Fritz Gesztesy,
Alexander Kiselev, and Konstantin A.~Makarov}
\address{Department of Mathematics,
University of
Missouri, Columbia, MO
65211, USA}
\email{fritz@math.missouri.edu\newline
\indent{\it URL:}
http://www.math.missouri.edu/people/fgesztesy.html}
\address{Department of Mathematics, University of
Missouri, Columbia, MO
65211, USA}
\email{makarov@math.missouri.edu\newline
\indent{\it URL:}
http://www.math.missouri.edu/people/kmakarov.html}
\address{Department of Mathematics, University of Chicago, 
5734 South University Avenue,  
Chicago, IL 60637-1546 USA}
\email{kiselev@math.uchicago.edu}
\thanks{To appear in {\it Math. Nachr.}}

\begin{abstract}
Let $g(z,x)$ denote the diagonal Green's matrix of a self-adjoint 
$m\times m$ matrix-valued Schr\"odinger operator $H=
-\f{d^2}{dx^2}I_m +Q(x)$ in $L^2 (\bbR)^{m}$, $m\in\bbN$. One 
of the principal results proven in this paper states that for a 
fixed $x_0\in\bbR$ and all $z\in\bbC_+$, $g(z,x_0)$ and 
$g^\prime (z,x_0)$ uniquely determine the matrix-valued
$m\times m$ potential $Q(x)$ for a.e.~$x\in\bbR$. We also prove 
the following local version of this result. Let $g_j(z,x)$, 
$j=1,2$ be the diagonal Green's matrices of the self-adjoint
Schr\"odinger  operators $H_j=-\f{d^2}{dx^2}I_m +Q_j(x)$ in $L^2
(\bbR)^{m}$. Suppose that for fixed $a>0$
and $x_0\in\bbR$, $\|g_1(z,x_0)-g_2(z,x_0)\|_{\bbC^{m\times m}}+
\|g_1^\prime (z,x_0)-g_2^\prime (z,x_0)\|_{\bbC^{m\times m}} 
\underset{|z|\to\infty}{=}O\big(e^{-2\Im(z^{1/2})a}\big)$ for $z$
inside  a cone along the imaginary axis with vertex zero and
opening  angle less than $\pi/2$, excluding the real axis. Then 
$Q_1(x)=Q_2(x)$ for a.e.~$x\in [x_0-a,x_0+a]$. 

Analogous results are proved for matrix-valued Jacobi and 
Dirac-type operators. \\

\noindent This is a revised and updated version of a 
previously archived file.
\end{abstract}

\maketitle

\section{Introduction}\lb{s1}

While various aspects of inverse spectral theory for scalar 
Schr\"odinger, Jacobi, and Dirac-type operators, and more
generally, for $2\times 2$ Hamiltonian systems, are 
well-understood by now, the corresponding theory for such 
operators and Hamiltonian systems with $m\times m$, $m\in\bbN$, 
matrix-valued coefficients is still in its infancy. A particular 
inverse spectral theory aspect we have in mind is that of
determining isospectral sets (manifolds) of such systems. It may, 
perhaps, come as a surprise that determining the isospectral set 
of Hamiltonian systems with matrix-valued periodic coefficients 
is a completely open problem. It appears to be no exaggeration to 
claim that (unless one considers trivial cases such as diagonal 
coefficient matrices, etc.) absolutely nothing seems to be known 
about the corresponding isospectral sets of periodic matrix-valued
Schr\"odinger operators with the sole exception of its 
compactness. The same ignorance applies to Jacobi, Dirac, and more
generally, to  periodic $2m\times 2m$ Hamiltonian systems with
$m\geq 2$. While the present paper is not sufficiently ambitious 
to change this sorry state of affairs, we will take a modest step
toward a closer investigation of inverse spectral problems and
prove a few uniqueness theorems for such systems, that is, 
determine spectral data (interpreted in a very broad sense) that
uniquely determine the matrix-valued coefficients in Schr\"odinger,
Jacobi, and Dirac-type systems. It should be mentioned that these 
types of problems are not just of interest in a spectral theoretic
context, but due to their implications for other areas such as
completely integrable systems (e.g., the nonabelian
Korteweg-deVries, Toda lattice, and Ablowitz-Kaup-Newell-Segur
(nonlinear Schr\"odinger) hierarchies), are also of considerable
interest to a much larger audience. 

Before we continue along this line of thought, it seems
appropriate to briefly mention some of the present day
knowledge of $2m\times 2m$ matrix-valued Hamiltonian systems and
the Weyl-Titchmarsh and spectral theory associated with them. To
save space we will often simultaneously discuss references on
Schr\"odinger, Dirac, and general Hamiltonian systems together
without differentiating between them, and occasionally single out
the finite-difference (Jacobi) systems. Moreover, there exists a
considerable amount of literature on the foundations of Hamiltonian
systems (and their special cases, such as Schr\"odinger, Jacobi,
and Dirac-type systems) which necessarily forces us to be rather
selective here. In particular, we focus primarily on the case
$m\geq 2$ for the remainder of this introduction. 

The basic Weyl-Titchmarsh theory of regular Hamiltonian systems can
be found in Atkinson's monograph \cite{At64}; Weyl-Titchmarsh
theory of singular Hamiltonian systems and their basic spectral
theory was developed by Hinton and Shaw and many others (see, e.g., 
\cite{Al95}, \cite{Hi69}--\cite{HSc97}, \cite{JNO00}, \cite{KR74},
\cite{Kr89a}, \cite{Kr89b}, \cite{LM00a}, \cite{Or76}, \cite{Ro60},
\cite{Ro69}, \cite{Sa92}, \cite{Sa94a}, \cite[Ch.~9]{Sa99a},
\cite{We87} and the references therein); the corresponding theory
for Jacobi systems can be found in \cite{Be68}, \cite{Fu76},
\cite{Sa97} and the literature therein. Various aspects of direct
spectral theory, including  investigations of the nature of the
spectrum involved, (regularized) trace  formulas, etc., appeared in 
\cite{BL99}, \cite{Ca98}, \cite{Ca99}, \cite{Ca00}, \cite{Cl90}, 
\cite{Cl94}, \cite{CGHL98}, \cite{CH98}, \cite{GH97}, 
\cite{GJ89}, \cite{Kr83}, \cite{LW99}, \cite{Ma65}, \cite{Pa95},
\cite{RH77}.  General asymptotic expansions of Weyl-Titchmarsh
matrices as the  (complex) spectral parameter tends to infinity
under optimal regularity assumptions on the coefficients are of
relatively recent origin and can be found in \cite{Cl92},
\cite{CG99}, \cite{CG01} (see also \cite{Sa88}, \cite{Tr00}). The
inverse scattering formalism for continuous Hamiltonian systems has
been studied by a variety of authors and we refer, for instance, to
\cite{AM63}, \cite{AG95}, \cite{AG98}, \cite{AK88}, \cite{AK90},
\cite{Ga68}, \cite{MO82}, \cite{MST01}, \cite{NJ55}, \cite{WK74},
\cite{Zh95}.  General inverse spectral theory, the existence of
transformation operators, etc., is discussed in \cite{GKS97},
\cite{LM00}, \cite{Ma94}, \cite{Ma99}, \cite{Sa90}, \cite{Sa97b}, 
\cite{Sa01}, \cite{Sa94a}, \cite{Sa96}, \cite{Sa99}, \cite{Sa99a}, 
and the references therein. Inverse monodromy problems for canonical
systems received a lot of attention recently. The interested reader
is referred to \cite{AD97}, \cite{AD00}, \cite{AD00a}, 
\cite{Ma95}--\cite{Ma99a}, \cite{Sa94}, \cite{Sa99a} and the
extensive literature cited therein. The corresponding inverse
spectral and scattering theory for matrix-valued finite difference
systems and its intimate connection to matrix-valued orthogonal
polynomials is treated in \cite{AG94}, \cite{AN84}, \cite{Be68},
\cite{DGK78}, \cite{DV95}--\cite{Fu76}, \cite{Ge81}, \cite{Ge82},
\cite{Lo99}, \cite{MBO92}, \cite{Sa97}.  More specific inverse
spectral problems, such as compactness of the isospectral set of
periodic Schr\"odinger operators \cite{Ca98a},  special isospectral
matrix-valued Schr\"odinger operators, and  Borg-type uniqueness
theorems (for periodic coefficients as well as eigenvalue problems
on compact intervals) were recently studied in \cite{Ch99a},
\cite{CS97}, \cite{CG01}, \cite{CGHL98}, \cite{De95}, \cite{JL98},
\cite{JL99}, \cite{Ma94}, \cite{Ma99}, \cite{Ma99a}, \cite{Sa99a},
\cite{SS98}. Moreover, direct spectral theory in the particular
case of periodic Hamiltonian systems (i.e., Floquet theory and
alike) has been studied in \cite{Ca98a}, \cite{Ca99}, \cite{CGHL98},
\cite{De80}, \cite{De95}, \cite{GL87}, \cite{Jo87}, \cite{Kr83a},
\cite{Kr83}, \cite{Ro63}, \cite{Ya92}--\cite{YS75b}, with many more
pertinent references to be found therein. 

Apart from Floquet theoretic applications in connection with
Hamiltonian systems already briefly touched upon, we also
need to mention applications to random Schr\"odinger operators
associated with strips as discussed, for instance, in \cite{CL90},
\cite{KS88}, \cite{La84}, \cite{La86}, and especially to nonabelian
completely integrable systems. Since the literature associated with
the latter topic is of enormous proportions, we can only refer to a
few pertinent publications, such as, \cite{AK90}, \cite{BGS86},
\cite{BG90}, \cite{Ch96}, \cite{Di91}--\cite{Du83}, \cite{GD77},
\cite{Ma78}, \cite{Ma88}, \cite{MO82}, \cite{OMG81}, \cite{Sa97a},
\cite{Sa88}--\cite{Sa94a}, \cite{Sa97}. The interested reader will 
find a wealth of additional material in these references.

Finally, we turn to the principal subject of this paper, that is,
uniqueness-type theorems for Schr\"odinger, Jacobi, and Dirac-type
operators. We have already mentioned a few uniqueness results in
connection with Borg-type theorems for such systems. Additional 
uniqueness-type results in terms of matrix-valued Wronskians and
transformation operators can be found in some papers by Leibenzon 
\cite{Le62}, \cite{Le62a} and Malamud \cite{Ma95}, \cite{Ma99}. 
The uniqueness theorems proven in this paper are directly 
formulated in terms of diagonal Green's matrices
$g(z,x_0)$ and their $x$-derivatives $g^\prime (z,x_0)$ at some
fixed $x_0\in\bbR$, for matrix-valued Schr\"odinger and Dirac-type
operators on $\bbR$ and similarly in terms of diagonal Green's
matrices $g(z,k_0,k_0)$ and off-diagonal Green's matrices
$G(z,k_0,k_0+1)$, etc., for some fixed $k_0\in\bbZ$, in the context
of matrix-valued Jacobi operators on $\bbZ$. Moreover, we prove
certain localized versions of these uniqueness theorems for
exponentially close diagonal Green's matrices with respect to $z$
as $|z|\to\infty$ in the Schr\"odinger and Dirac-type context and
analogous theorems in the case of Jacobi operators whose (diagonal
and certain off-diagonal) Green's matrices differ by an inverse
power of $z$ as
$|z|\to\infty$. To be a bit more specific, we briefly describe
some of our principal results in the case of self-adjoint
Schr\"odinger operators $H$ in $L^2(\bbR)^{m}$, $m\in\bbN$
associated with $m\times m$ matrix-valued differential expressions
$-\f{d^2}{dx^2}I_m +Q(x)$, with $Q=Q^*\in L^1_{\loc}
(\bbR)^{m\times m}$. 

Let $\bbC_+$ be the open complex upper
half-plane. Denoting by $g(z,x)$, $z\in\bbC\backslash\bbR$, 
$x\in\bbR$ the diagonal Green's matrix associated with $H$ (i.e.,
the integral kernel of the resolvent of $H$ on the diagonal, 
$g(z,x)=(H-z)^{-1}(x,x)$) and by $g^\prime (z,x)$ its
$x$-derivative, we will prove the following
result (cf.~Theorem~\ref{t3.5}). 

\begin{theorem} \lb{t1.1}
Fix $x_0\in\bbR$. Then $g(z,x_0)$ and $g^\prime (z,x_0)$ for 
all $z\in\bbC_+$ uniquely determine the matrix-valued $m\times m$
potential $Q(x)$ for a.e.~$x\in\bbR$.
\end{theorem}

Moreover, let $H_j$, $j=1,2$ be two Schr\"odinger operators 
associated with the $m\times m$ matrix-valued potentials $Q_j(x)$ 
and $g_j(z,x)$ the corresponding diagonal Green's matrices of
$H_j$. Then we prove the following local version of
Theorem~\ref{t1.1} (cf.~Theorem~\ref{t3.5a}).

\begin{theorem} \lb{t1.2}
Let $a>0$ and $x_0\in\bbR$. If 
\begin{align}
&\|g_1(z,x_0)-g_2(z,x_0)\|_{\bbC^{m\times m}}+
\|g_1^\prime (z,x_0)-g_2^\prime (z,x_0)\|_{\bbC^{m\times m}} 
\no \\
&\underset{|z|\to\infty}{=}O(e^{-2\Im(z^{1/2})a}) \lb{1.1}
\end{align}
along a ray with $0<\arg(z)<\pi$, then
\begin{equation}
Q_1(x)=Q_2(x) \text{ for a.e.~$x\in [x_0-a,x_0+a]$}. \lb{1.2}
\end{equation}
\end{theorem}

Our results were of course inspired by analogous ones in the
special scalar context $m=1$, but also by a very interesting
uniqueness theorem proven by Berezanskii in 1953 \cite{Be53}, 
\cite{Be64} in the context of multi-dimensional Schr\"odinger
operators $H$. Berezanskii's point of departure is a bit different
from ours. He considered the spectral kernel
$\vartheta(\lambda,x,y)$ of $H$ (i.e., the integral 
kernel of the spectral projection $E_H (\lambda)$ of $H$, 
$\vartheta(\lambda,x,y)=E_H (\lambda,x,y)$) and 
appropriate normal derivatives of it with respect to $x$ and
$y$ across an arbitrarily small piece of a smooth surface in 
$\bbR^d$, $d=2,3$.
There are additional hypotheses in Berezanskii's work which 
need not be discussed here. In a certain sense we replaced the
spectral projection of $H$ by its resolvent and (in our
one-dimensional context) the arbitrarily small piece of a surface
by the point $x_0$ and the corresponding normal derivative by 
$d/dx$.

In Section~\ref{s2} we review the basic Weyl-Titchmarsh theory for
matrix-valued Schr\"    odinger and Dirac-type operators as needed
in our principal Section~\ref{s3}. The latter contains our
uniqueness results for Schr\"odinger and Dirac-type operators. The
final  Section~\ref{s4} then treats Weyl-Titchmarsh theory and
uniqueness theorems for matrix-valued Jacobi operators. 

\section{Matrix-Valued Schr\"{o}dinger and Dirac-Type 
Operators} \lb{s2}

In this section we briefly recall the Weyl-Titchmarsh 
theory
for matrix-valued Schr\"{o}dinger  and Dirac-type 
operators.  In 
order to
treat both cases in
parallel, we use the fact that both are special cases 
of Hamiltonian
systems and
hence develop the theory from that point of view. Throughout
this paper all matrices will be considered over the field of
complex numbers $\bbC$ and the corresponding linear space 
of $k\times\ell$ matrices will be denoted by 
$\bbC^{k\times\ell}$.

The basic assumption for Sections~\ref{s2} and \ref{s3} of 
this paper will be the following.

\begin{hypothesis}\lb{h2.1}
Fix $m\in\bbN$ and define the $2m\times 2m$ matrix
\begin{equation}
J=\begin{pmatrix}0& -I_m \\ I_m & 0  \end{pmatrix}. \lb{2.1}
\end{equation}
We consider two principal cases: Either \\
$(i)$ suppose $Q=Q^*\in L_{\loc}^1(\bbR)^{m\times m}$ and
introduce the $2m\times 2m$ matrices
\begin{equation}
A=\begin{pmatrix} I_m & 0 \\ 0& 0\end{pmatrix}, \quad
B(x)=\begin{pmatrix}-Q(x) & 0 \\0 & I_m\end{pmatrix}, 
\lb{2.2}
\end{equation}
or \\
$(ii)$  suppose
\begin{equation}
A=I_{2m}, \quad B = B^*\in L_{\loc}^1(\bbR)^{2m\times 2m}. 
\lb{2.3}
\end{equation}
\end{hypothesis}

Given Hypothesis~\ref{h2.1} we consider the Hamiltonian 
system
\begin{equation}
J\varPsi^\prime(z,x)=(zA+B(x))\varPsi(z,x)  
\text{ for a.e. $x\in\bbR$},  \lb{2.4}
\end{equation}
where $z\in\bbC$ plays the role 
of a spectral
parameter and $\Psi(z,x)$ is
assumed to satisfy
\begin{equation}
\varPsi(z,\dott) \in \AC_{\loc}(\bbR)^{2m\times 2m}. \lb{2.5}
\end{equation}
Here and later on, $I_p$ denotes the identity matrix in
$\bbC^{p\times p}$ for 
$p\in\bbN$, $M^*$
the adjoint (i.e., complex
conjugate transpose), $M^t$ the transpose of the matrix 
$M$, and
$\AC_{\loc}(\bbR)$ denotes the
set of locally absolutely continuous functions on $\bbR$.  
At times it will be convenient to
consider a $2m\times r$ solution matrix of \eqref{2.4}, with
$r=1,\dots,2m$, whose entries will then be assumed to lie in 
$\AC_{\loc}(\bbR)$.

One verifies that Hypothesis~\ref{h2.1}\,(i) governs the 
case of matrix-valued Schr\"{o}dinger
operators. In fact, inserting 
\begin{equation}
\varPsi(z,x)=\begin{pmatrix}\psi_1(z,x) & 0 \\
\psi_2(z,x) & 0 \end{pmatrix} 
\lb{2.8}
\end{equation}
into equation \eqref{2.4} then yields
\begin{align}
-\psi_1^{\prime\prime}(z,x)+Q(x)\psi_1(z,x)& 
=z\psi_1(z,x), \lb{2.6} \\
\psi_2(z,x)&=\psi_1^\prime(z,x). \lb{2.7}
\end{align}
Here it is assumed that
\begin{equation}
\psi_j(z,\dott)\in\AC_{\loc}(\bbR)^{m\times m}, \quad j=1,2.
\lb{2.9}
\end{equation}
In the case of Hypothesis~\ref{h2.1}\,(ii), \eqref{2.4} 
represents a Dirac-type system.

In connection with our uniqueness result for Dirac-type 
operators in Section~\ref{s3} we will make the 
following additional assumption due to the freedom of certain 
gauge transformations in the Dirac context, which leaves the
corresponding spectral matrix, but not the potential $B(x)$
invariant, as discussed by Gasymov \cite{Ga68} and Gasymov and
Levitan \cite{GL66}.  

\begin{hypothesis} \lb{h2.2}
Suppose Hypothesis~\ref{h2.1} and in the case of Dirac-type
operators {\rm (}i.e., assuming  \eqref{2.3}{\rm )} suppose 
that $B$ is locally essentially bounded, $B\in
L^\infty_{\loc}(\bbR)^{2m\times 2m}$, and that
$B(x)$  is of the special form
\begin{align}
& B(x)=\begin{pmatrix}
B_{1,1}(x)& B_{1,2}(x) \\ B_{1,2}(x) & -B_{1,1}(x)
\end{pmatrix}, \lb{2.13a} \\
& \text{ with }
B_{1,1}(x)=B_{1,1}(x)^*, \,\, B_{1,2}(x)=B_{1,2}(x)^*, 
\,\, j=1,2. \no
\end{align}
\end{hypothesis}
Equation \eqref{2.13a} represents a typical normal 
form inspired by the requirement
\begin{equation}
JB(x)+B(x)J=0 \text{ for a.e. } x\in\bbR. \lb{2.13b}
\end{equation}
Various possible normal forms for Dirac-type operators are 
discussed in \cite{GL66}, \cite{LM00}, \cite{Ma94}, \cite{Ma99},
\cite{Ma65}, and in the monographs \cite[Ch.~9]{LS75},
\cite[Ch.~7]{LS91}, \cite[p.~193--195]{Ma86}. The continuity 
assumption on $B$ in Hypothesis~\ref{h2.2} is further 
discussed in the paragraph following Theorem~\ref{t2.5}.

Next we briefly turn to Weyl-Titchmarsh theory associated 
with \eqref{2.4} and recall some
of the results developed by Hinton and Shaw in a series 
of papers devoted to spectral theory of
(singular) Hamiltonian systems \cite{HS81}--\cite{HS86}
(see also \cite{Kr89a},
\cite{Kr89b}). While they discuss
\eqref{2.4} under much more general hypotheses on $A(x)$ 
and $B(x)$, we confine ourselves here to
the special cases of  matrix-valued Schr\"{o}dinger and 
Dirac-type systems governed by Hypothesis~\ref{h2.1}. The following 
facts on Weyl-Titchmarsh theory are taken from \cite{CG99} and
\cite{CG01} and hence we omit the corresponding proofs.

 Let $\Psi(z,x,x_0)$ be a normalized fundamental system of 
solutions of \eqref{2.4} at some
$x_0\in\bbR$, that is, $\Psi(z,x,x_0)$ satisfies
\begin{equation}
J\Psi^\prime(z,x)=(zA+B(x))\Psi(z,x), \quad z\in\bbC 
\lb{2.14}
\end{equation}
for a.e.\ $x\in\bbR$, and
\begin{equation}
\Psi(z,x_0,x_0)=I_{2m}. \lb{2.15}
\end{equation}
Moreover, we partition $\Psi(z,x,x_0)$ as
\begin{align}
\Psi(z,x,x_0)&=\big(\psi_{j,k}(z,x,x_0)\big)_{1\leq j,k\leq 2} 
\no \\
&=\begin{pmatrix} \Theta(z,x,x_0) & \Phi(z,x,x_0)\end{pmatrix}
=\begin{pmatrix}\theta_1(z,x,x_0) &
\phi_1(z,x,x_0)\\
\theta_2(z,x,x_0)& \phi_2(z,x,x_0)\end{pmatrix}, \lb{2.16}
\end{align}
where $\theta_j(z,x,x_0)$ and $\phi_j(z,x,x_0)$ for $j=1,2$ 
are $m\times m$ matrices, entire with
respect to $z\in\bbC$, and normalized according to 
\eqref{2.15}, that is,
\begin{equation}
\theta_1(z,x_0,x_0)=\phi_2(z,x_0,x_0)=I_m, \quad
\phi_1(z,x_0,x_0)=\theta_2(z,x_0,x_0)=0. \lb{2.17}
\end{equation}
Next, let $\beta_{j} \in \bbC^{m\times m}$, $j=1,2$, 
introduce the matrix 
$\beta=(\beta_{1}\,\beta_{2})\in\bbC^{2m\times m}$, and assume that
\begin{equation}
\rank(\beta)=m, \quad \beta\beta^*=I_m, \,
 \text{ and either } \Im(\beta_2\beta_1^*)\leq 0 \text{ or } 
\Im(\beta_2\beta_1^*)\geq 0. \lb{2.18}
\end{equation}

\noindent One can then prove the following result.

\begin{lemma} [see, e.g.,~\cite{CG01}] \lb{l2.4} Assume
Hypothesis~\ref{h2.1}, let $\Theta(z,x,x_0)$ and
$\Phi(z,x,x_0)$ be defined as in \eqref{2.16}, and suppose 
$\beta$ satisfies \eqref{2.18}. Then, for
$c\in\bbR\backslash\{x_0\}$, $\beta\Phi(z,c,x_0)$ is singular if and
only if $z$ is an eigenvalue for the  regular  boundary value
problem given by \eqref{2.4}, \eqref{2.5} together with the
separated boundary conditions
\begin{equation}
(I_m\, 0)\varPsi(z,x_0)=0, \quad \beta\varPsi(z,c)=0. \lb{BC}
\end{equation}
\end{lemma}

\noindent For the regular boundary value problem described
in Lemma~\ref{l2.4}, the boundary conditions in \eqref{BC} are
self-adjoint whenever $\Im(\beta_2\beta_1^*)=0$.

Lemma~\ref{l2.4} provides appropriate conditions for defining  a
certain meromorphic $m\times m$ matrix $M(z,c,x_0,\beta)$.

\begin{definition}\lb{dMF}
Assume Hypothesis~\ref{h2.1} and let $ \Theta(z,x,x_0)$, and
$\Phi(z,x,x_0)$ be defined as in \eqref{2.16} with $\beta$
satisfying \eqref{2.18}. For
$c\ne x_0$, and $\beta\Phi(z,c,x_0)$ nonsingular let
\begin{equation} 
M(z,c,x_0,\beta) =
-[\beta\Phi(z,c,x_0)]^{-1}[\beta\Theta(z,c,x_0)]. \lb{MF}
\end{equation}
$M(z,c,x_0,\beta)$ is said to be the
{\it Weyl-Titchmarsh $M$-function} for the regular boundary value
problem described in Lemma~\ref{l2.4}.
\end{definition}
The Weyl-Titchmarsh $M$-function in \eqref{MF} is an
$m\times m$ matrix-valued function with meromorphic entries whose
poles correspond to eigenvalues for the regular boundary value
problem given by \eqref{2.4}, \eqref{2.5}, and \eqref{BC}.
Moreover, if $M\in \bbC^{m\times m}$, and one defines
\begin{equation}
U(z,x,x_0)= \begin{pmatrix}
u_1(z,x,x_0)\\u_2(z,x,x_0) \end{pmatrix}=
\Psi(z,x,x_0)\begin{pmatrix}I_m\\M\end{pmatrix}, \lb{2.20}
\end{equation}
with $u_j(z,x,x_0)\in \bbC^{m\times m}$, $j=1,2$, then
$U(z,x,x_0)$ will satisfy the boundary condition at $x=c$ in
\eqref{BC} whenever $M=M(z,c,x_0,\beta)$. Intimately connected 
with the matrices introduced in Definition~\ref{dMF} is the set of
$m\times m$ complex matrices known as the Weyl disk. 

To describe this set, we first introduce the matrix-valued function 
$E(M)$: With $c\ne x_0$, $z\in\bbC\backslash\bbR$, and with
$U(z,c,x_0)$  defined by \eqref{2.14} in terms of a  matrix
$M\in\bbC^{m\times m}$, let
\begin{equation} \lb{2.380}
E(M) = \sigma(x_0,c,z)U(z,c,x_0)^*(iJ)U(z,c,x_0),
\end{equation}
where 
\begin{equation}
\sigma(s,t,z)=\frac{(s-t)\Im (z)}{|(s-t)\Im (z)|},\quad
\sigma(s,t)=\sigma(s,t,i),\quad \sigma(z)= \sigma(1,0,z), \lb{2.19}
\end{equation}
with $s,t\in\bbR$, $s\ne t$. 

\begin{definition}\lb{dWD}
Assume Hypothesis~\ref{h2.1} and fix $x_0\in\bbR$,
$c\in\bbR\backslash\{x_0\}$, and $z\in\bbC\backslash\bbR$. Then 
$D(z,c,x_0)$ will denote the collection of all $M\in
\bbC^{m\times m}$ for which $E(M)\le 0$, where $E(M)$ is
defined in
\eqref{MF}. $D(z,c,x_0)$ is said to be a {\it Weyl disk}.
The set of $M\in \bbC^{m\times m}$ for which $E(M) = 0$ is
said to be a {\it Weyl circle} (even when $m>1$).
\end{definition}

This definition leads to a  presentation  
that is a generalization of the description first given 
by Weyl~\cite{We10}; a presentation which is geometric in nature, 
involves the contractive matrices
$V\in\bbC^{m\times m}$, such that $VV^*\le I_m$, and provides the
justification for the geometric terms of circle and disk and 
(see, e.g., \cite{HS81}, \cite{HSc93}, \cite{Kr89a}, \cite{Or76}).

The disk has also been characterized 
in terms of matrices which statisfy hypothesis \eqref{2.18} and
which  serve as boundary data for the regular boundary value
problem  described in Lemma~\ref{l2.4} (see, e.g. \cite{At88a},
\cite{At88}). In particular, $\calD(z,c,x_0)$ denotes the
collection of all $M\in\bbC^{m\times m}$ obtained by the
construction given in \eqref{MF} where $c\ne x_0$, $z\in\bbC/\bbR$,
where $\beta$ is  the $m\times m$ matrices defined in hypothesis
\eqref{2.18} for which $\sigma(c,x_0,z)\Im(\beta_2\beta_1^*)\ge 0$.

We note that the Weyl circle corresponds to the regular
boundary value problems in Lemma~\ref{l2.4} with separated, 
self-adjoint boundary conditions. 

\begin{lemma}[\cite{HS84}, \cite{HSc93}, \cite{Kr89a}]\lb{l2.11}
Let $M\in\bbC^{m\times m}$, $c\ne x_0$, and $z\in\bbC\backslash\bbR$.
Then, $E(M)=0$ if and only if there is a $\beta\in \bbC^{m\times
2m}$ satisfying $\beta J\beta^*=0$ and $\beta\beta^*=I_m$ such that
\begin{equation}\lb{2.24}
0=\beta U(z,c,x_0),
\end{equation}
where $U(z,c,x_0)$ is defined in \eqref{2.14} in
terms of $M$.  With $\beta$ so defined,
\begin{equation}\lb{2.25}
M=-[\beta\Phi(z,c,x_0)]^{-1}[\beta\
\Theta(z,c,x_0)],
\end{equation}
that is, $M=M(z,c,x_0,\beta)$.
\end{lemma}

Next, we recall a fundamental property associated with matrices
in $\calD(z,c,x_0)$.

\begin{lemma} [see, e.g., \cite{CG01}]  \lb{l2.8} 
If $M\in\calD(z,c,x_0)$, then 
\begin{equation}\lb{Hgz}
\sigma(c,x_0,z) \Im (M)>0.
\end{equation}
Moreover, whenever $\beta\in\bbC^{m\times 2m}$ satisfies 
$\beta J\beta^*=0$ and $\beta\beta^*=I_m$,
\begin{equation}\lb{2.270}
M(\bar z,c,x_0,\beta)= M(z,c,x_0,\beta)^*.
\end{equation}
\end{lemma}

For $c>x_0$, the function $M(z,c,x_0,\beta)$, 
defined by \eqref{MF}, and satisfying \eqref{Hgz}, is said 
to be a matrix-valued {\it Herglotz} function of rank $m$. Hence,
for $\Im(\beta_2\beta_1^*)=0$,  poles of
$M(z,c,x_0,\beta)$, $c>x_0$, are at most of first order, are
real, and have nonpositive residues. Such functions admit
representations of the form
\begin{equation}\lb{NP}
M(z,c,x_0,\beta)=C_1 + zC_2 + \int_{-\infty}^\infty
d\Omega(\lambda,c,x_0,\beta)\,\left(
\frac{1}{\lambda-z} -\frac{\lambda}{1+\lambda^2} \right), 
\quad c>x_0, 
\end{equation}
where (with $\beta$ fixed) 
$C_2\ge 0$ and $C_1$ are constant $m\times m$ self-adjoint
matrices, and $\Omega(\lambda,c,x_0,\beta)$ is a nondecreasing
$m\times m$ matrix-valued function such that
\begin{align}
&\int_{-\infty}^{\infty}||
d\Omega(\lambda,c,x_0,\beta)||_{\bbC^{m\times m}} \,
(1+\lambda^2)^{-1} <\infty, \lb{NPa} \\
&\Omega((\lambda, \mu],c,x_0,\beta)= \lim_{\delta\downarrow 0}
\lim_{\epsilon \downarrow 0}\frac{1}{\pi}\int_{\lambda
+ \delta}^{\mu
+ \delta }d\nu\,  \sigma(c,x_0,\nu +i\epsilon)\Im\left( M(\nu
+i\epsilon,c,x_0,\beta)\right). \lb{NPb}
\end{align}
For the self-adjoint boundary value problems in connection
with Schr\"odinger, Dirac, and Jacobi operators discussed in this
paper, $C_2=0$ in \eqref{NP} (and also later in \eqref{2.25a} and
\eqref{2.42}) and $\Omega(\lambda,c,x_0,\beta)$ is then piecewise
constant with jump discontinuities at the  eigenvalues of the
problem. Analogous statements apply to $-M(z,c,x_0,\beta)$ if
$c<x_0$.

We  further note that the sets $D(z,c,x_0)$ are closed, 
convex, and are nested with respect to increasing or decreasing 
values of  $c$ (cf. \cite{HS81}, \cite{HS83}, \cite{HS84},
\cite{Kr89a}, \cite{Or76}), that is,
\begin{equation}
D(z,c_2,x_0)\subseteq D(z,c_1,x_0) \quad
\text{for}\quad x_0<c_1\le c_2\quad \text{or}
\quad c_2\le c_1<x_0. \lb{2.20b}
\end{equation}
Hence, the intersection of this nested sequence, as $c\to \pm
\infty$, is nonempty, closed and convex.  We say that this
intersection is a limiting set for the nested sequence.
\begin{definition}\lb{dLWD}
Assume Hypothesis~\ref{h2.1} and let  $D_\pm(z,x_0)$ denote the
closed, convex set in the space of $m\times m$ matrices which is
the limit, as
$c\to \pm\infty$,  of the nested collection of sets
$D(z,c,x_0)$ given in Definition~\ref{dWD}.
$D_\pm(z,x_0)$ is said to be a limiting {\em disk}.
Elements of $D_\pm(z,x_0)$
are denoted by $M_\pm(z,x_0)\in \bbC^{m\times m}$.
\end{definition}

In light of the containment described in
\eqref{2.20b},  for $c\ne x_0$ and $z\in \bbC\backslash\bbR$,
\begin{equation}
D_\pm(z,x_0)\subseteq D(z,c,x_0). \lb{2.20c}
\end{equation}
Given its characterization in the second paragraph following 
Definition~\ref{dWD}, we see that $D_\pm(z,x_0)$ consists
of those matrices $M(z,c,x_0,\beta)$ where
$\beta$ satisfies \eqref{2.18}. In particular,
\begin{equation}
M_\pm(z,x_0)=M(z,c,x_0,\beta) \lb{2.20d}
\end{equation}
for an appropriate choice of $\beta$.

When $D_\pm(z,x_0)$ is a singleton matrix, the
system \eqref{2.4} is said to be in the {\it limit point} (l.p.) 
case at
$\pm\infty$.  When $D_\pm(z,x_0)$ has nonempty interior, then
\eqref{2.4} is said to be in the {\it limit circle} (l.c.) 
case at
$\pm\infty$. Indeed, for the case $m=1$, the limit circle case
corresponds to $D_\pm(z,x_0)$ being a disk in $\bbC$.

These apparent geometric properties for the disk
correspond to  analytic properties for the solutions of the 
Hamiltonian system \eqref{2.4}, \eqref{2.5}.  To recall this
correspondence, we introduce the following spaces in which we assume
that $ -\infty\le a< b \le \infty$,
\begin{subequations}\lb{2.20e}
\begin{align}
L_A^2((a,b))&=\bigg\{\phi:(a,b)\to\bbC^{2m}  \bigg| \int_a^b dx\,
(\phi(x),A\phi(x))_{\bbC^{2m}}<\infty \bigg\}, \lb{2.20ea}
\\
N(z,\infty)&=\{\phi\in L_A^2((c,\infty)) \mid J\phi^\prime
=(zA+B)\phi
\text{ a.e. on $(c,\infty)$} \}, \lb{2.20eb}
\\
N(z,-\infty)&=\{\phi\in L_A^2((-\infty,c)) \mid
J\phi^\prime=(zA+B)\phi
\text{ a.e. on $(-\infty,c)$} \}, \lb{2.20ec}
\end{align}
\end{subequations}
for some $c\in\bbR$ and $z\in\bbC$. (Here
$(\phi,\psi)_{\bbC^n}=\sum_{j=1}^n \overline\phi_j\psi_j$
denotes the standard scalar product in $\bbC^n$, abbreviating
$\chi\in\bbC^n$ by
$\chi=(\chi_1,\dots,\chi_n)^t$.)  Both dimensions of the
spaces in \eqref{2.20eb} and \eqref{2.20ec},
$\dim_\bbC(N(z,\infty))$ and $\dim_\bbC(N(z,-\infty))$,
are constant for $z\in\bbC_\pm=\{\zeta\in\bbC
\mid \pm\Im(\zeta)> 0 \}$ (see, e.g., \cite{At64},
\cite{KR74}). One then observes that  the Hamiltonian
system \eqref{2.4} is in the limit point case at $\pm\infty$
whenever
\begin{equation}
\dim_\bbC(N(z,\pm\infty))=m \text{  for all
$z\in\bbC\backslash\bbR$} \lb{2.20f}
\end{equation}
and in the limit circle case at $\pm\infty$  whenever 
\begin{equation}
\dim_\bbC(N(z,\pm\infty))=2m \text{  for all $z\in\bbC$.} \lb{2.20g}
\end{equation}

Next we recall the fact that the Dirac-type systems considered in
this paper  are always in the limit point case at $\pm\infty$.

\begin{lemma} \lb{l2.15}
Assume Hypothesis~\ref{h2.1}. Then the limit point case holds 
for Dirac-type systems {\rm (}i.e., for $A(x)=I_{2m}$ in
\eqref{2.4}{\rm )} at $\pm \infty$.
\end{lemma}
Lemma~\ref{l2.15}, under varying sets of assumptions on $B(x)$, is
well-known  to experts in the field. For instance, in the case
$m=1$ and with $B_{1,2}(x)=B_{2,1}(x)$ this fact can be found in
\cite{We71}. For $B\in C(\bbR)^{2m\times 2m}$ and a more general
constant matrix $A$,  this result is proven in \cite{LM00} (their
proof, however,  extends to the current $B\in L^1_{\loc} (\bbR)$
case). More generally, multi-dimensional  Dirac operators with
$L^2_{\loc} (\bbR^n)$-type coefficients (and  additional
conditions) can be found in \cite{LO82}. A short proof of
Lemma~\ref{l2.15} can be found in \cite{CG01}. After completion 
of this paper we became aware of a recent
preprint by Lesch and Malamud \cite{LM00a} which provides a thorough
study of self-adjointness questions for more general Hamiltonian 
systems than those considered in this paper.

In either the limit point or limit circle cases,
$M_\pm(z,x_0)\in \partial D_{\pm}(z,x_0)$ is said to be a
{\em half-line Weyl-Titchmarsh  matrix}. (In the l.p.~case one of
course has $D_{\pm}(z,x_0)=\partial D_{\pm}(z,x_0)$.) Each such
matrix is  associated with the construction of a self-adjoint
operator acting on
$L_A^2([x_0,\pm \infty))\cap \AC([x_0,\pm\infty))$ for the 
Hamiltonian
system \eqref{2.4}, \eqref{2.5}. However, for those intermediate
cases where $m<\dim_\bbC(N(z,\pm\infty))<2m$, Hinton and Schneider
have shown that not every element of the boundary $\partial
D_{\pm}(z,x_0)$ of $D_{\pm}(z,x_0)$ is a half-line Weyl-Titchmarsh
matrix, and have characterized those elements of the boundary that
are (cf. \cite{HSc93}, \cite{HSc97}).

For later reference we summarize the principal results on
$M_{\pm}(z,x_0)$ in the following theorem.

\begin{theorem}
[\cite{GT97}, \cite{HS81}, 
\cite{HS82}, \cite{HS86}, \cite{KS88}] \lb{thm2.12} 
Assume Hypothesis~\ref{h2.1} and suppose \, that 
$z\in\bbC\backslash\bbR$, $x_0\in\bbR$. Then \\
$(i)$ $\pm M_{\pm}(z,x_0)$ is an $m\times m$ matrix-valued 
Herglotz function of maximal rank. In particular,
\begin{gather}
\Im(\pm M_{\pm}(z,x_0)) > 0, \quad z\in\bbC_+, 
\lb{2.21} \\
M_{\pm}(\overline z,x_0)=M_{\pm}( z,x_0)^*,\lb{2.22} \\
\rank (M_{\pm}(z,x_0))=m, \lb{2.23} \\
\lim_{\varepsilon\downarrow 0} M_{\pm}(\lambda+
i\varepsilon,x_0) \text{
exists for a.e.\
$\lambda\in\bbR$}. \lb{2.24a}
\end{gather}
Local singularities of $\pm M_{\pm}(z,x_0)$ and 
$\mp M_{\pm}(z,x_0)^{-1}$ are necessarily real and at most of first
order in the sense that 
\begin{align}
&\mp \lim_{\epsilon\downarrow0}
\left(i\epsilon\,
M_{\pm}(\lambda+i\epsilon,x_0)\right) \geq 0, \quad \lambda\in\bbR,
\lb{2.24b} \\ 
& \pm \lim_{\epsilon\downarrow0}
\left(\f{i\epsilon}{M_{\pm}(\lambda+i\epsilon,x_0)}\right)
\geq 0, \quad \lambda\in\bbR. \lb{2.24c}
\end{align}
$(ii)$  $\pm M_{\pm}(z,x_0)$ admit the representation
\begin{equation}
\pm M_{\pm}(z,x_0)=F_\pm(x_0)+\int_\bbR 
d\Omega_\pm(\lambda,x_0) \,
\big((\lambda-z)^{-1}-\lambda(1+\lambda^2)^{-1}\big), \lb{2.25a} 
\end{equation}
where
\begin{equation}
F_\pm(x_0)=F_\pm(x_0)^*, \quad \int_\bbR 
\|d\Omega_\pm(\lambda,x_0)\|_{\bbC^{m\times m}}\,
(1+\lambda^2)^{-1}<\infty \lb{2.27} 
\end{equation}
and
\begin{equation}
\Omega_\pm((\lambda,\mu],x_0)=\lim_{\delta\downarrow
0}\lim_{\varepsilon\downarrow 0}\f1\pi
\int_{\lambda+\delta}^{\mu+\delta} d\nu \, \Im(\pm
M_\pm(\nu+i\varepsilon,x_0)). \lb{2.29} 
\end{equation}
$(iii)$  Define the $2m\times m$ matrices
\begin{align}
\Psi_\pm(z,x,x_0)&=\begin{pmatrix}\psi_{\pm,1}(z,x,x_0)\\
\psi_{\pm,2}(z,x,x_0)  \end{pmatrix} \no \\
&=\begin{pmatrix}\theta_1(z,x,x_0) & \phi_1(z,x,x_0)\\
\theta_2(z,x,x_0)& \phi_2(z,x,x_0)\end{pmatrix} 
\begin{pmatrix} I_m \\
M_\pm(z,x_0) \end{pmatrix}, \lb{2.31}
\end{align}
then 
\begin{equation}
\Im(M_\pm(z,x_0))=\Im(z) \int_{x_0}^{\pm\infty}dx\,
\Psi_\pm(z,x,x_0)^* A
\Psi_\pm(z,x,x_0), \lb{2.32}
\end{equation}
where $A$ has been introduced in \eqref{2.2}, \eqref{2.3}.
\end{theorem}

In order to describe the Green's matrix associated with 
\eqref{2.14}, we assume the hypotheses of  
Theorem~\ref{thm2.12} and introduce
\begin{align}
K(z,x,x^\prime)=\Psi_\mp(z,x,x_0)[M_-(z,x_0) &
-M_+(z,x_0)]^{-1}\Psi_\pm(\overline z,x^\prime,x_0)^*, 
\no \\
& \hspace*{2.05cm} x\lessgtr x^\prime,\, 
z\in\bbC\backslash\bbR, \lb{2.33}
\end{align}
and
\begin{align}
M(z,x_0)&=\big(M_{j,j^\prime}(z,x_0)\big)_{j,j^\prime=1,2} 
\no \\
&=[K(z,x_0,x_0+0)+K(z,x_0,x_0-0)]/2,
\quad  z\in\bbC\backslash\bbR,  \lb{2.34} \\
M_{1,1}(z,x_0)&=[M_-(z,x_0)-M_+(z,x_0)]^{-1}, \no \\
M_{1,2}(z,x_0)&=2^{-1} [M_-(z,x_0)-M_+(z,x_0)]^{-1}
[M_-(z,x_0)+M_+(z,x_0)], \no \\
M_{2,1}(z,x_0)&=2^{-1}
[M_-(z,x_0)+M_+(z,x_0)][M_-(z,x_0)-M_+(z,x_0)]^{-1},\no \\
M_{2,2}(z,x_0)&=M_\pm(z,x_0)[M_-(z,x_0)-M_+(z,x_0)]^{-1}
M_\mp(z,x_0). \no 
\end{align}
Next let $\phi\in L_A^2(\bbR)$ and consider
\begin{equation}
J\psi^\prime(z,x)=(zA+B(x))\psi(z,x)+A \phi(x), 
\quad z\in\bbC\backslash\bbR
\lb{2.36}
\end{equation}
for a.e.\ $x\in\bbR$.  Then \eqref{2.36} has a 
unique solution
$\psi(z,\dott)\in
L_A^2(\bbR)\cap\AC_{\loc}(\bbR)^{2m}$ given by \cite{HS81}, 
\cite{HS83}
\begin{equation}
\psi(z,x)=\int_\bbR dx^\prime\, K(z,x,x^\prime) 
A\phi(x^\prime). \lb{2.37}
\end{equation}
In the following we associate an operator $T$ in 
$L^2(\bbR)^m$,
respectively, $L^2(\bbR)^{2m}$ with
the Hamiltonian system \eqref{2.4} denoted by 
\begin{equation}
T=\begin{cases}H & \text{if $A=\begin{pmatrix} 
I_m & 0 \\0 & 0
\end{pmatrix}$}, \\
D & \text{if $A=I_{2m}$}. \end{cases} \lb{2.38}
\end{equation}
Here $H$ is defined by 
\begin{equation}
((H-z)^{-1}f)(x)= \int_\bbR dx^\prime\, 
K_{1,1}(z,x,x^\prime)f(x^\prime), 
\quad z\in\bbC\backslash\bbR, \, f\in
L^2(\bbR)^m, \lb{2.38a} 
\end{equation}
and this results in the following matrix-valued Schr\"{o}dinger 
operator in
$L^2(\bbR)^m$,
\begin{align}
H&=-I_m \f{d^2}{dx^2}+Q, \lb{2.39} \\
\dom(H)&=\{g\in L^2(\bbR)^m \mid g,g^\prime\in
\AC_{\loc}(\bbR)^m;\, \text{s. s.-a. b. c. at\,}\pm\infty;
\no \\
& \hspace*{5.1cm} (-I_m g^{\prime\prime}+Qg)\in L^2(\bbR)^m\}, \no
\end{align}
where $\text{``s. s.-a. b. c. at }\pm\infty$'' denotes a 
separated self-adjoint boundary condition at $+\infty$ 
and/or $-\infty$ 
(if any) induced by $M_+(z,x_0)$ and/or $M_-(z,x_0)$. Here
$K_{1,1}(z,x,x^\prime)$  is the left upper $m\times m$ 
submatrix of $K(z,x,x^\prime)$, obtained from the partitioning
$K(z,x,x^\prime)=(K_{j,k}(z,x,x^\prime))_{j,k=1}^2$.
Similarly, 
$D$ defined by 
\begin{equation}
((D-z)^{-1}\psi)(x)= \int_\bbR dx^\prime\, 
K(z,x,x^\prime)\psi(x^\prime), 
\quad z\in\bbC\backslash\bbR, \, \psi\in
L^2(\bbR)^{2m}, \lb{2.39a} 
\end{equation}
represents the following Dirac-type operator 
in $L^2(\bbR)^{2m}$,
\begin{align}
D&=J \f{d}{dx}-B, \lb{2.40} \\
\dom(D)&=\{\phi\in L^2(\bbR)^{2m}\mid \phi
\in\AC_{\loc}(\bbR)^{2m}; \,(J\phi^\prime-B\phi)\in 
L^2(\bbR)^{2m} \}, \no
\end{align}
where we took into account the limit point property of Dirac-type
operators as described in Lemma~\ref{l2.15}. Thus, $H$ and $D$ are
self-adjoint operators in
$ L^2(\bbR)^m$  and $ L^2(\bbR)^{2m}$, respectively.

The basic results on $M(z,x_0)$ then read as follows.

\begin{theorem} [\cite{GT97}, \cite{HS81}, \cite{HS82}, 
\cite{HS86}, \cite{KS88}] \lb{thm2.4} 
Assume Hypothesis~\ref{h2.1} and suppose \, that $z\in\bbC 
\backslash \bbR$, $x_0\in\bbR$.  Then \\
$(i)$ $M(z,x_0)$ is a matrix-valued Herglotz function of rank 
$2m$ with
representation
\begin{equation}
M(z,x_0)=F(x_0)+\int_\bbR d\Omega(\lambda,x_0)\,
\big((\lambda-z)^{-1}-\lambda(1+\lambda^2)^{-1}\big), \lb{2.42} 
\end{equation}
where
\begin{equation}
F(x_0)=F(x_0)^*, \quad \int_\bbR \Vert 
d\Omega(\lambda,x_0)
\Vert_{\bbC^{2m\times 2m}} \,(1+\lambda^2)^{-1}<\infty \lb{2.44}
\end{equation}
and
\begin{equation}
\Omega((\lambda,\mu],x_0)=\lim_{\delta\downarrow
0}\lim_{\varepsilon\downarrow 0}\f1\pi
\int_{\lambda+\delta}^{\mu+\delta} d\nu \, 
\Im(M(\nu+i\varepsilon,x_0)). \lb{2.46}
\end{equation}
$(ii)$ $z\in\rho(T)$ if and only if $M(z,x_0)$ is holomorphic 
near $z$. 
\end{theorem}

The fundamental uniqueness theorem for $T$ in terms of
$M(z,x_0)$ then reads as follows.

\begin{theorem} [\cite{CG99}, \cite{CG01}, \cite{GS99},
\cite{LM00}, \cite{Ro60}] \lb{t2.5} Let $x_0\in\bbR$ and assume 
Hypothesis~\ref{h2.2}. Then 
$M(z,x_0)$ for all $z\in\bbC_+$ uniquely determines T. In
particular, in the case of the Schr{\"o}\-dinger operator $H$,
$M(z,x_0)$, $z\in\bbC_+$, uniquely determines $Q(x)$ for 
a.e.~$x\in\bbR$,  and in the case of the Dirac-type operator $D$,
$M(z,x_0)$, $z\in\bbC_+$, uniquely determines
$B(x)=\left(\begin{smallmatrix}B_{1,1}(x)& B_{1,2}(x)
\\B_{1,2}(x) & -B_{1,1}(x)\end{smallmatrix}\right)$ 
for a.e.~$x\in\bbR$.
\end{theorem}

In the context of Schr\"odinger operators, Theorem~\ref{t2.5}
represents the matrix-valued extension of the celebrated
Borg-Marchenko  uniqueness theorem \cite{Bo52}, \cite{Ma50},
\cite{Ma52}  (see also \cite{GS00}, \cite{Si98}) for Schr\"odinger
operators on a  half-line in the scalar case 
$m=1$.~The matrix Schr\"odinger case has been treated by 
Rofe-Beketov \cite{Ro60} (assuming $Q$ to be continuous on 
$\bbR$) and more recently in \cite{GS99} (in the general 
case $Q\in L^1([0,R])^{m\times m}$ for all $R>0$).  In the
Dirac case, and for $B$ continuous, Theorem~\ref{t2.5} is 
contained in \cite{LM00}. The Dirac result under the weaker
Hypothesis~\ref{h2.2} is due to \cite{CG01}. In the scalar case
$m=1$, boundedness of $B$ is immaterial for the uniqueness result
to hold in  Theorem~\ref{t2.5} (cf.~the discussion in \cite{CG01}). 
However, for $m\geq 2$, the corresponding  uniqueness result appears
to have been proved only assuming 
$B\in L^\infty_{\loc}(\bbR)^{2m\times 2m}$ thus far. This is the
reason why we added  the local essentially boundedness assumption of
$B$ to Hypothesis~\ref{h2.2}.  As soon as one succeeds in removing
this $L^\infty_{\loc}$-assumption of $B$ in Theorem~\ref{t2.5}, it
can be removed everywhere in this  paper. 

The recent new results on local Borg-Marchenko uniqueness
theorems for matrix-valued Schr\"odinger and Dirac-type 
operators on a half-line in \cite{CG01}, \cite{GS99} 
then yield the following local extension of Theorem~\ref{t2.5}.

\begin{theorem} [\cite{CG01}, \cite{GS99}] \lb{t2.6}
Suppose $B_j$, $j=1,2$, satisfy Hypothesis~\ref{h2.2} and
let $a>0$, $a_\pm >0$, and $x_0\in\bbR$. \\
$(i)$ Denote by $M_{j,\pm}(z,x_0)$, $j=1,2$ the 
$m\times m$ Weyl-Titchmarsh matrices associated with the
corresponding Hamiltonian systems \eqref{2.4} on 
$(x_0,\pm\infty)$. In the Schr\"odinger case assume that 
\begin{equation}
\|M_{1,\pm}(z,x_0)-M_{2,\pm}(z,x_0)\|_{\bbC^{m\times
m}}\underset{|z|\to\infty}{=}
O\big(e^{-2\Im(z^{1/2})a_\pm}\big)
 \lb{2.48}
\end{equation}
along a ray with $0<\arg(z)<\pi$. In the Dirac-type
case suppose that 
\begin{equation}
\|M_{1,\pm}(z,x_0)-M_{2,\pm}(z,x_0)\|_{\bbC^{m\times
m}}\underset{|z|\to\infty}{=}
O\big(e^{-2\Im(z)a_\pm}\big) \lb{2.48a}
\end{equation}
along a ray with $0<\arg(z)<\pi/2$ and along a ray with
$\pi/2<\arg(z)<\pi$. Then, in either case, 
\begin{equation}
B_1(x)=B_2(x) \text{ for a.e. } x\in [x_0,x_0\pm a_{\pm}]. 
\lb{2.49}
\end{equation}
$(ii)$ Denote by $M_j(z,x_0)$, $j=1,2$ the 
$2m\times 2m$ Weyl-Titchmarsh matrices associated with the
corresponding Hamiltonian systems \eqref{2.4} on $\bbR$. 
In the Schr\"odinger case assume that 
\begin{equation}
\|M_1(z,x_0)-M_2(z,x_0)\|_{\bbC^{2m\times
2m}}\underset{|z|\to\infty}{=}
O\big(e^{-2\Im(z^{1/2})a}\big) \lb{2.50}
\end{equation}
along a ray with $0<\arg(z)<\pi$. In the Dirac-type
case suppose that
\begin{equation}
\|M_1(z,x_0)-M_2(z,x_0)\|_{\bbC^{2m\times
2m}}\underset{|z|\to\infty}{=}
O\big(e^{-2\Im(z^{1/2})a}\big) \lb{2.50a}
\end{equation}
along a ray with $0<\arg(z)<\pi/2$ and along a ray with
$\pi/2<\arg(z)<\pi$. Then, in either case, 
\begin{equation}
B_1(x)=B_2(x) \text{ for a.e. } x\in [x_0-a,x_0+a]. 
\lb{2.51}
\end{equation}
\end{theorem}
\begin{proof}
In the context of Schr\"odinger operators, case\,(i) is proved in
\cite{GS99}; in the corresponding Dirac-type situation case\,(i) is 
proved in \cite{CG01}. Case\,(ii) then follows by  combining
\eqref{2.34} and case\,(i), taking into account  the asymptotic
expansions (along any ray with $0<\arg(z)<\pi$)
\begin{equation}
M_\pm (z,x_0)\underset{|z|\to\infty}{=}\pm iz^{1/2}I_m 
+ o(1) \lb{2.52}
\end{equation}
in the case of Schr\"odinger operators (cf.~\cite{CG99}) and
\begin{equation}
M_\pm (z,x_0)\underset{|z|\to\infty}{=}\pm iI_m + o(1) \lb{2.53}
\end{equation} 
in the case of Dirac-type operators (cf.~\cite{CG01}).
\end{proof}

In the case of scalar Schr\"odinger operators, 
Theorem~\ref{t2.6}\,(i) is due to Simon \cite{Si98}. An 
alternative proof was presented in \cite{GS99}, and more recently
in \cite{Be00}. Most recently, Alexander Sakhnovich kindly informed
us that his integral representation of the Weyl-Titchmarsh matrix in
\cite{Sa88a} can be used to derive asymptotic expansions for the
Weyl-Titchmarsh matrix and its associated matrix-valued spectral
function, and also yields a result analogous to
Theorem~\ref{t2.6}\,(i) for a certain class of canonical systems.
Moreover, in the case of skew-adjoint Dirac-type systems, similar
results are discussed in
\cite{Sa90} and applied to the nonlinear Schr\"odinger equation on
a half-axis. 

\section{Uniqueness Results for Schr\"odinger and 
Dirac-Type Operators} \lb{s3}

In this section we prove our uniqueness results for
matrix-valued Schr\"odinger and Dirac-type operators.

Since in the case of Schr\"odinger operators, $B(x)$ 
is of a very special nature,
\begin{equation}
B(x)=\begin{pmatrix}-Q(x)& 0
\\ 0 & I_m \end{pmatrix}, \quad Q(x)^*=Q(x), \lb{3.0}  
\end{equation}
and only the $m\times m$ submatrix 
$B_{1,1}=-Q(x)$ contains information on $Q(x)$,
we will focus on the $m\times m$ submatrix
$M_{1,1}(z,x_0)$  of $M(z,x_0)$ in
\eqref{2.34} and \eqref{2.42}, \eqref{2.44}. By 
\eqref{2.33} and
\eqref{2.38a} one infers that the
Green's matrix $G(z,x,x^\prime)$ of $H$ is given by
\begin{align}
G(z,x,x^\prime)&=K_{1,1}(z,x,x^\prime) \no \\
&=\psi_{\mp,1}(z,x,x_0)[M_-(z,x_0)-M_+(z,x_0)]^{-1}
\psi_{\pm,1}(\overline z,x^\prime,x_0)^*, \no \\
& \hspace*{6.3cm} x\lesseqqgtr x^\prime, \, 
z\in\bbC\backslash\bbR, \lb{3.1}
\end{align}
where the $m\times m$ matrices $\psi_{\pm,j}(z,x,x_0)$ 
are defined in
\eqref{2.31}, that is,
\begin{equation}
\psi_{\pm,j}(z,x,x_0)=\theta_j(z,x,x_0)+
\phi_j(z,x,x_0)M_\pm(z,x_0), \quad
j=1,2. \lb{3.2}
\end{equation}
For future purpose we abbreviate the diagonal Green's function 
$G(z,x,x)$ of $H$ by $g(z,x)$, that is,
\begin{equation}
g(z,x)=G(z,x,x), \quad z\in\bbC\backslash\bbR, \,\,
x\in\bbR. 
\lb{3.2a}
\end{equation}
Moreover, since ``diagonal elements'' of matrix-valued 
Herglotz functions
are (lower-dimensional) matrix-valued Herglotz
functions  (cf., e.g., \cite{GT97}),
$M_{1,1}(z,x_0)$ $=g(z,x_0)$
is an $m\times m$ matrix-valued Herglotz function 
satisfying
\begin{align}
M_{1,1}(z,x_0)&=g(z,x_0)=K_{1,1}(z,x_0,x_0) \lb{3.3} \\
&=[M_-(z,x_0)-M_+(z,x_0)]^{-1} \lb{3.3a} \\
&= F_{1,1}(x_0)+\int_\bbR d\Omega_{1,1}(\lambda,x_0)\,
\big((\lambda-z)^{-1}-\lambda(1+\lambda^2)^{-1}\big), \lb{3.4} 
\end{align}
where
\begin{equation}
F_{1,1}(x_0)=F_{1,1}(x_0)^*, \quad \int_\bbR \Vert
d\Omega_{1,1}(\lambda,x_0)
\Vert_{\bbC^{m\times m}} \, (1+\lambda^2)^{-1}<\infty \lb{3.6}  
\end{equation}
and
\begin{equation}
\Omega_{1,1}((\lambda,\mu],x_0)=\lim_{\delta\downarrow 0}
\lim_{\varepsilon\downarrow
0}\f1\varepsilon \int_{\lambda+\delta}^{\mu+\delta}d\nu\, 
\Im(M_{1,1}(\nu+i\varepsilon,x_0)). \lb{3.8} 
\end{equation}

In complete analogy to the Schr\"odinger case, we also 
introduce
\begin{equation}
g(z,x_0)=M_{1,1}(z,x_0)=[M_-(z,x_0)-M_+(z,x_0)]^{-1}, \quad 
z\in\bbC\backslash\bbR \lb{3.10}
\end{equation}
in the case of Dirac-type operators $D$ and recall that  
$B(x)$ is now of the type 
\begin{align}
&B(x)=\begin{pmatrix}
B_{1,1}(x)& B_{1,2}(x) \\ B_{2,1}(x) & B_{2,2}(x)
\end{pmatrix}, \lb{3.11} \\
& \text{ with }
B_{1,2}(x)^*=B_{2,1}(x), \,\, B_{j,j}(x)^*=B_{j,j}(x),
\,\, j=1,2. \no
\end{align}
Next, we provide a quick derivation of
the matrix Riccati-type equations satisfied by $M_\pm(z,x)$ in
the Schr\"odinger and Dirac cases. This result is not new
(see, e.g., \cite{CGHL98}, \cite{GH97}, \cite{Jo87}). We
only mention its proof in some detail to reach a certain level
of completeness and since it is a fundamental ingredient in
our principal uniqueness result, Theorem~\ref{t3.5}. 

\begin{lemma}\lb{l3.1}
Assume Hypothesis~\ref{h2.1}, $z\in\bbC\backslash\bbR$, and
$x\in\bbR$.  Then in the case of matrix-valued Schr\"odinger
operators $H$, $M_\pm(z,x)$ satisfies the standard
Riccati-type equation,
\begin{equation}
M_\pm^\prime(z,x)+M_\pm(z,x)^2=Q(x)-z I_m \lb{3.17}
\end{equation}
for a.e.~$x\in\bbR$. In the case of Dirac-type operators $D$,
$M_\pm(z,x)$ satisfies
\begin{align}
&M_\pm^\prime(z,x)+zM_\pm(z,x)^2+M_\pm(z,x)B_{2,2}(x)M_\pm(z,x)
+B_{1,2}(x)M_\pm(z,x)  \no \\
&+M_\pm(z,x)B_{2,1}(x)=-B_{1,1}(x)-z I_m \lb{3.18}
\end{align}
for a.e.~$x\in\bbR$.
\end{lemma}
\begin{proof}
In order to prove \eqref{3.17}, let
\begin{equation}
\hat \Psi(z,x,x_0)=\begin{pmatrix} \psi_{-,1}(z,x,x_0)
&\psi_{+,1}(z,x,x_0)\\
\psi_{-,2}(z,x,x_0)& \psi_{+,2}(z,x,x_0)\end{pmatrix}, \quad
z\in\bbC\backslash\bbR \lb{3.12}
\end{equation}
be the fundamental system of solutions of \eqref{2.14} 
as defined in
\eqref{2.31} and
observe that
\begin{equation}
\psi_{+,2}(z,x,x_0)=\psi_{+,1}^\prime(z,x,x_0) \lb{3.13}
\end{equation}
in the Schr\"{o}dinger operator case. Thus, any 
nonnormalized solutions
$\ti \psi_{\pm,1}(z,\dott)\in L^2((\pm\infty,c))^{m\times m}$, 
$c\in\bbR$ of
$-\psi_1^{\prime\prime}+Q\psi_1=z\psi_1$ for $z\in\bbC
\backslash\bbR$ are of the
type
\begin{equation}
\ti \psi_{\pm,1}(z,x)=\psi_{\pm,1}(z,x,x_0)C_\pm \lb{3.14}
\end{equation}
for some nonsingular $m\times m$ matrices $C_\pm$.  Thus,
\begin{equation}
\ti \psi_{\pm,1}(z,x_0)=C_\pm, \quad \ti
\psi_{\pm,1}^\prime(z,x_0)=M_\pm(z,x_0)C_\pm \lb{3.15}
\end{equation}
by \eqref{2.17} and \eqref{2.31}.  In particular,
\begin{equation}
M_\pm(z,x_0)=\ti \psi_{\pm,1}^\prime(z,x_0) 
\ti \psi_{\pm,1}(z,x_0)^{-1}, \quad
z\in\bbC\backslash\bbR
\lb{3.16}
\end{equation}
is independent of the normalization chosen for 
$\ti \psi_{\pm,1}(z,x_0)$.
Varying the reference
point
$x_0\in\bbR$ then yields \eqref{3.17}. Similarly, in order to
derive \eqref{3.18}, let
\begin{equation}
\hat \Psi(z,x,x_0)=\begin{pmatrix} \psi_{-,1}(z,x,x_0)
&\psi_{+,1}(z,x,x_0)\\
\psi_{-,2}(z,x,x_0)& \psi_{+,2}(z,x,x_0)\end{pmatrix}, 
\quad z\in\bbC_+
\lb{3.29}
\end{equation}
be the fundamental system of solutions of \eqref{2.14} 
as defined in
\eqref{2.31}.
Any other fundamental system of \eqref{2.14} of the form
\begin{equation}
\tilde \Psi(z,x)=\begin{pmatrix} \ti \psi_{-,1}(z,x)&
\ti \psi_{+,1}(z,x)\\
\ti \psi_{-,2}(z,x)& \ti \psi_{+,2}(z,x)\end{pmatrix}, 
\quad z\in\bbC_+,
\lb{3.30}
\end{equation}
with $\tilde \Psi_{\pm}(z,x)=(\ti \psi_{\pm,1}(z,x),\ti
\psi_{\pm,2}(z,x))^t$ a basis of
$N(z,\pm\infty)$ for
$z\in\bbC_+$, necessarily must be of the form
\begin{equation}
\tilde \Psi(z,x)= \hat \Psi(z,x,x_0)C. \lb{3.31}
\end{equation}
Here $C$ is of the type
\begin{equation}
C=\begin{pmatrix}C_+& 0 \\ 0 & C_-  \end{pmatrix}  \lb{3.32}
\end{equation}
and $C_\pm$ are invertible $m\times m$ matrices.  In 
particular,
\begin{equation}
\ti \psi_{-,j}(z,x)=\psi_{-,j}(z,x,x_0)C_-, \quad
\ti \psi_{+,j}(z,x)=\psi_{+,j}(z,x,x_0)C_+, \,\, j=1,2
\lb{3.33}
\end{equation}
and thus,
\begin{equation}
\ti \psi_{\pm,1}(z,x_0)=C_\pm, \quad \ti \psi_{\pm,2}(z,x_0)
=M_\pm(z,x_0)C_\pm
\lb{3.34}
\end{equation}
by \eqref{2.17} and \eqref{2.31}.  Consequently,
\begin{equation}
M_\pm(z,x_0)=\ti \psi_{\pm,2}(z,x_0) 
\ti \psi_{\pm,1}(z,x_0)^{-1}, \quad
z\in\bbC_+\lb{3.35}
\end{equation}
is independent of the normalization chosen for 
$\ti \psi_{\pm,j}$.  Varying
$x_0\in\bbR$ one
verifies \eqref{3.18}.
\end{proof}
\begin{corollary} \lb{c3.2}
Assume Hypothesis~\ref{h2.1}, $z\in\bbC\backslash\bbR$, and
$x\in\bbR$. Moreover, define 
\begin{equation}
g(z,x)=[M_-(z,x)-M_+(z,x)]^{-1}, \quad z\in\bbC\backslash\bbR, 
\,\, x\in\bbR \lb{3.37}
\end{equation} 
as in \eqref{3.2a} and \eqref{3.10}. Then in the case of
matrix-valued Schr\"odinger operators, the diagonal Green's matrix
$g(z,x)$ of $H$, satisfies
\begin{align}
g'(z,x)&=g(z,x)M_+(z,x)+M_-(z,x)g(z,x),
\lb{3.38} \\
g'(z,x) \mp
I_m&=g(z,x)M_\pm(z,x)+M_\pm(z,x)g(z,x) \lb{3.39}
\end{align}
for all $x\in\bbR$. In the case of Dirac-type operators $D$,
$g(z,x)$ satisfies
\begin{align}
g'(z,x)&=zg(z,x)M_+(z,x)+zM_-(z,x)g(z,x)
-g(z,x)M_\pm(z,x)B_{1,1}(x) \no \\
&\quad -B_{1,1}(x)M_\mp(z,x)g(z,x)+
g(z,x)B_{1,2}(x)+B_{1,2}(x)g(z,x) \lb{3.40}
\end{align}
for a.e.~$x\in\bbR$.
\end{corollary}
\begin{proof}
In order to prove \eqref{3.38} one only needs to observe
\begin{align}
&(d/dx)\big(g(z,x)^{-1}\big)=
-g(z,x)^{-1}g'(z,x)\
g(z,x)^{-1}\no \\
&=M_-'(z,x)-M_+'(z,x) \no \\
&=M_+(z,x)^2-M_-(z,x)^2 \no \\
&=M_+(z,x)[M_+(z,x)-M_-(z,x)]  \no \\
&\quad \, +[M_+(z,x)-M_-(z,x)]M_-(z,x) \no \\
&=-M_+(z,x)g(z,x)^{-1}-g(z,x)^{-1}M_-(z,x) \lb{3.41}
\end{align}
for a.e.~$x\in\bbR$, using \eqref{3.17} and \eqref{3.37}. By
continuity of the right-hand side of \eqref{3.41} with
respect to $x$, \eqref{3.41} extends to all $x\in\bbR$. Equation
\eqref{3.41}  then proves \eqref{3.38}, and using \eqref{3.37}
again, 
\eqref{3.38} proves \eqref{3.39}. Similarly, to prove
\eqref{3.40}, one computes
\begin{align}
g'(z,x)&=-[M_-(z,x)-M_+(z,x)]^{-1}[M_-'(z,x)-M_+'(z,x)]
\times \no \\
& \quad \times [M_-(z,x)-M_+(z,x)]^{-1} \no \\
& =-[M_-(z,x)-M_+(z,x)]^{-1}[zM_+(z,x)^2-zM_-(z,x)^2  \no \\ 
& \quad -M_+(z,x)B_{1,1}(x)M_+(z,x)+M_-(z,x)B_{1,1}(x)M_-(z,x) 
\no \\
& \quad +B_{1,2}(x)M_+(z,x) -B_{1,2}(x)M_-(z,x)
+M_+(z,x)B_{1,2}(x) \no \\
& \quad -M_-(z,x)B_{1,2}(x)][M_-(z,x)-M_+(z,x)]^{-1}
\no \\
& =zg(z,x)M_+(z,x)+zM_-(z,x)g(z,x) \no \\
& \quad +[M_-(z,x)-M_+(z,x)]^{-1}[M_+(z,x)B_{1,1}(x)M_+(z,x)
\no \\
& \quad +M_\pm(z,x)B_{1,1}(x)M_\mp(z,x)
-M_\pm(z,x)B_{1,1}(x)M_\mp(z,x) \no \\
& \quad +M_-(z,x)B_{1,1}(x)M_-(z,x)][M_-(z,x)-M_+(z,x)]^{-1}
\no \\
& \quad +g(z,x)B_{1,2}(x)+B_{1,2}(x)g(z,x) \no \\
& =zg(z,x)M_+(z,x)+zM_-(z,x)g(z,x)
-g(z,x)M_\pm(z,x)B_{1,1}(x) \no \\
&\quad -B_{1,1}(x)M_\mp(z,x)g(z,x)+
g(z,x)B_{1,2}(x)+B_{1,2}(x)g(z,x) \lb{3.42}
\end{align}
for a.e.~$x\in\bbR$, using \eqref{3.18} and \eqref{3.37}.
\end{proof}

Denoting by $\spec(A)$ the spectrum (i.e., the set of
eigenvalues) of an $n \times n$ matrix 
$A\in\bbC^{n\times n}$, we
next recall a well-known lemma concerning the solution of
Sylvester's matrix equation. We also recall that $A$ is 
said to be accretive if 
$\Re((\ul x,A\ul x)_{\bbC^n})\geq 0$ for all 
$\ul x\in\bbC^n$. Finally, we recall that $|||\cdot |||$ 
is said to be a unitarily invariant norm on 
$\bbC^{n\times n}$ if $|||UAV|||=|||A|||$ for all 
$U,V\in\bbC^{n\times n}$ unitary and all 
$A\in\bbC^{n\times n}$. 

\begin{lemma} [c.f.~\cite{Bh97},
Theorems~VII.2.1 and  VII.2.4 and \cite{BDM83}] \lb{l3.3}   
Assume $A,B,C \in \bbC^{n\times n}$ for some $n\in\bbN$. \\
$(i)$ Suppose that
\begin{equation}
\spec(A)\cap\spec(B)=\emptyset. \lb{3.43}
\end{equation}
Then the equation
\begin{equation}
AX-XB=C \lb{3.44}
\end{equation}
has a unique solution $X\in\bbC^{n \times n}$. Moreover, 
suppose 
$\Gamma$ is a closed counterclockwise oriented Jordan 
contour such 
that $\spec (A)$ has winding number $+1$ and $\spec (B)$ 
has winding number $0$ with respect to $\Gamma$. Then an 
explicit solution for $X$ in \eqref{3.44} is provided by 
\begin{equation}
X=(2\pi i)^{-1}\oint_\Gamma d\zeta\,
(A-\zeta)^{-1}C(B-\zeta)^{-1}. 
\lb{3.44a}
\end{equation}
$(ii)$ Assume that $A-(\delta/2)I_n$ and $-B-(\delta/2)I_n$ 
are accretive. Then \eqref{3.44} has a unique solution
$X\in\bbC^{n\times n}$ satisfying 
\begin{equation}
|||X|||\leq \delta^{-1}|||C||| \lb{3.44aa}
\end{equation}
for any unitarily invariant norm $|||\cdot |||$ on 
$\bbC^{n\times n}$. 
\end{lemma}

Without going into details, we note that Lemma~\ref{l3.3}\,(i) 
remains valid for bounded operators on a Banach space and 
also extends to certain cases where $A$ and $B$ are unbounded 
(cf., e.g., \cite{Ph91}). Similarly, Lemma~\ref{l3.3}\,(ii) is 
proven for maximal accretive operators $A-(\delta/2)I$, 
$-B-(\delta/2)I$ in a Hilbert space in \cite{BDM83}. 

Next, we recall one more auxiliary result to prove 
the principal uniqueness results of this section,
Theorems~\ref{t3.5} and \ref{t3.5a}. Although we are only 
applying the following result in the special case of $m\times m$ 
matrices, we state the corresponding facts in the 
infinite-dimensional context since, in contrast to 
Lemma~\ref{l3.3}, this appears to be a less familiar result. In the
following let $\calH$ be a complex separable Hilbert space with
scalar product $(\cdot,\cdot)_\calH$ (and norm $\|\cdot \|_\calH$),
and  denote by $\calB(\calH)$ the Banach space of bounded linear 
operators on $\calH$.

\begin{lemma} [c.f. \cite{GMN99}, Lemma~2.3]
\lb{l3.4} Let $M(z)$ be an operator-valued Herglotz
function with values in $\calB(\calH)$ {\rm(}i.e., $M$ is 
analytic in $\bbC_+$ with $\Im(M(z))\geq 0${\rm)}. 
If $(\Im(M(z_0)))^{-1}\in\calB(\calH) $ for some
$z_0\in\bbC_+$,  then $(\Im(M(z)))^{-1}\in\calB(\calH) $ 
for all $z\in\bbC_+$ and 
\begin{equation}
\spec(M(z))\subset\bbC_+ \text{ and hence } \,
\spec(M(z))\cap\bbR =\emptyset. \lb{3.45a}
\end{equation}
\end{lemma}
\begin{proof}
Suppose there is a $z_1\in \bbC_+$ such that 
$(\Im(M(z_1)))^{-1}\notin\calB(\calH)$. Then there is 
a sequence $\{e_n\}_n\subset \calH$, 
$\|e_n\|_\calH=1$, $n\in \bbN$, such that
\begin{equation}\lb{3.45b}
\lim_{n\to \infty}\|\Im(M(z_1))e_n\|_\calH=0.
\end{equation}
Any Herglotz  function $M(z)$ with values in $\calB(\calH)$ 
admits the representation
\begin{equation}\lb{3.45c}
M(z)=A+Bz+R^{1/2} (I_{\calK}+zL)(L-z)^{-1}R^{1/2}\big\vert_{\calH}, 
\end{equation}
where $A=A^*\in \calB(\calH)$, $0\le B\in \calB(\calH)$, $\calK$ is
a Hilbert space, $\calK\subseteq \calH$, $L=L^*$ is a self-adjoint
operator in $\calK$, and $0\leq R\in\calB(\calK)$ with
$R\vert_{\calK \ominus \calH}=0$. By \eqref{3.45c},
\begin{align}
\Im (M(z))&=\Im (z)( B+R^{1/2}(I_{\calK}+L^2)
\big((L-\Re(z))^2+(\Im (z))^2)^{-1}R^{1/2}\big\vert_{\calH}\big) \no
\\ &= \Im (z) \big(B+R^{1/2}C(z)R^{1/2}\big\vert_{\calH}\big),
\label{3.45d}
\end{align}
where $C(z)=(I_{\calK}+L^2)((L-\Re(z))^2+(\Im (z))^2)^{-1}>0$
is invertible in $\calK$, $C(z)^{-1}\in\calB(\calK)$, $ z\in 
\bbC_+$.
 Thus, \eqref{3.45b} together with \eqref{3.45d} imply 
\begin{equation}
\lim_{n\to \infty}
\big((e_n, Be_n)_{\calH} + (R^{1/2}\big\vert_{\calH} e_n, C(z_1)
R^{1/2}\big\vert_{\calH} e_n)_{\calH}\big)= 0 \lb{3.45e}
\end{equation}
and hence $\lim_{n\to \infty} \|B^{1/2}e_n\|_{\calH} = 
\lim_{n\to \infty}\|Be_n\|_{\calH}= 
\lim_{n\to \infty}$ $\|R^{1/2}\vert_{\calH}e_n\|_{\calH}= 0$.
Applying \eqref{3.45d} again one infers
$\lim_{n\to \infty}\|\Im(M(z))e_n\|_{\calH}= 0  \text{ for all }
z\in \bbC_+$, 
contradicting the hypothesis $(\Im(M(z_0)))^{-1}\in\calB(\calH)$.

Since $\Im(M(z))$, $z\in \bbC_+$ is invertible, the numerical range
$W_z$  of the operator $M(z)$
\begin{equation}
W_z=\{ (M(z)e,e)_\calH \, | \, e\in \calH, \|e\|_\calH=1\}
\lb{3.45f}
\end{equation}
is a subset of the half-plane
$\{ \zeta\in \bbC_+ \, | \, 
\Im(\zeta)\ge \|(\Im(M(z))^{-1})\|_{\calB(\calH)}^{-1} \}$. By the
spectral inclusion theorem  (see, e.g., \cite[Theorem1.2-1]{GR97}),
the spectrum of $M(z)$ is contained in  the closure  of the
numerical range $W_z$ and hence  \eqref{3.45a} holds.
\end{proof} 

Moreover, we will also apply the following elementary result. 

\begin{lemma}\lb{l.ric}
Let $A, B \in \bbC^{n\times n}$ for some $n\in \bbN$ and
\begin{equation}\lb{AB}
\|A\|_{\bbC^{n\times n}}<1/2, \quad\|B\|_{\bbC^{n\times n}}<1/2
\end{equation} for some matrix norm $\|\cdot \|$ on 
$\bbC^{n\times n}$.  Then the matrix-valued Riccati equation
\begin{equation}\lb{j1}
Y=YAY+B    
\end{equation}
 has a unique solution  $X\in \bbC^{n\times n}\cap \ol K$, where 
$\ol K=\{X\in\bbC^{n\times n} \,|\, \|X\|_{\bbC^{n\times n}}\le 1\}
$ is the closed unit ball in $\bbC^{n\times n}$. 
Moreover, the solution $X$ is   given by the convergent
series 
\begin{equation}
  X=\sum_{k=0}^\infty B_k, \lb{j.2}
\end{equation}
where 
\begin{equation}
B_0=B, \quad B_{k+1}=B_kAB_k+B, \quad k\in \bbN_0. \lb{j.3}
\end{equation}
\end{lemma}
\begin{proof}  Introducing the  map 
\begin{equation}
F:\bbC^{n\times n}\longrightarrow \bbC^{n\times n}, \quad 
F(X)=XAX+B, 
\end{equation}
the  estimate 
\begin{equation}
\|F(X)\|_{\bbC^{n\times n}}=\|B+XAX\|_{\bbC^{n\times n}}\le
\|B\|_{\bbC^{n\times n}}+\|X\|_{\bbC^{n\times
n}}^2\|A\|_{\bbC^{n\times n}}<1, \quad X\in \ol K,
\end{equation}
combined with \eqref{AB} shows that   
$F$ maps $\ol K$ into itself. Moreover, using
\begin{align}
 &\|F(X)-F(Y)\|_{\bbC^{n\times n}}=\|XAX-YAY\|_{\bbC^{n\times n}} 
\no \\
&=\|XAX-XAY+XAY-YAY\|_{\bbC^{n\times n}} \no \\ 
& \le(\|X\|_{\bbC^{n\times n}}+\|Y\|_{\bbC^{n\times n}})\, 
\|A\|_{\bbC^{n\times n}}\|X-Y\|_{\bbC^{n\times n}} \no \\
&\le 2\|A\|_{\bbC^{n\times n}}\|X-Y\|_{\bbC^{n\times n}}, 
\quad X,Y\in \ol K,
\end{align}
\eqref{AB} implies that 
$F$ is a strict contraction. An application of  Banach's fixed point
theorem then proves the existence of a unique fixed point $X\in
\ol K$. Moreover, the fixed point $X$ can be obtained as the limit
of the sequence $\{B_k\}_{k\in\bbN}$, where 
\begin{equation}
B_{k+1}=F(B_k), \quad k\in\bbN_0, 
\end{equation}
and any fixed $B_0\in\ol K$, 
\begin{equation}
X=\lim_{k\to\infty}F(B_k).
\end{equation}
Taking $B_0=B$ proves representations \eqref{j.2} and \eqref{j.3},
which in turn completes the proof since any fixed point of $F$ is a
solution of the matrix-valued Riccati equation \eqref{j1}. 
\end{proof}

\begin{corollary}\lb{cor}
Assume $A_j, $ $B_j$, $j=1,2$, satisfy the assumptions  of 
Lemma~\ref{l.ric} and $X_j$, $j=1,2$, are the corresponding
(unique) solutions of the matrix-valued Riccati equations
\begin{equation}
Y=YA_jY+B_j, \quad j=1,2
\end{equation}
in the closed unit ball $\ol K$. Then
\begin{equation}\lb{j8}
\|X_1-X_2\|_{\bbC^{n\times n}}\le\frac{\|A_1
-A_2\|_{\bbC^{n\times n}}+\|B_1
-B_2\|_{\bbC^{n\times n}}}{1-2\min_{j=1,2}
\|A_j\|_{\bbC^{n\times n}}}.
\end{equation}
\end{corollary}
\begin{proof} For $X_j\in \ol K, $ $j=1,2$ the estimate 
\begin{align}
&\|X_1-X_2\|_{\bbC^{n\times n}}=\|B_1-B_2
+X_1A_1X_1-X_2A_2X_2\|_{\bbC^{n\times n}} \no \\
&\le \|B_1-B_2\|_{\bbC^{n\times n}}
+\|X_1A_1X_1-X_2A_2X_2\|_{\bbC^{n\times n}} \no \\
&=\|B_1-B_2\|_{\bbC^{n\times n}}
+\|X_1A_1X_1-X_2A_1X_2\|_{\bbC^{n\times n}} \no \\
& \quad +\|X_2A_1X_2-X_2A_2X_2\|_{\bbC^{n\times n}} \no \\
&\le \|B_1-B_2\|_{\bbC^{n\times n}}+2\|A_1\|_{\bbC^{n\times n}}
\|X_1-X_2\|_{\bbC^{n\times n}} \no \\
& \quad +\|X_2\|_{\bbC^{n\times n}}^2 
\|A_1-A_2\|_{\bbC^{n\times n}} \no \\
&\le\|A_1-A_2\|_{\bbC^{n\times n}}+\|B_1-B_2\|_{\bbC^{n\times n}}
+2\|A_1\|_{\bbC^{n\times n}}\|X_1-X_2\|_{\bbC^{n\times n}}
\end{align}
shows that
\begin{equation}\lb{j9}
\|X_1-X_2\|_{\bbC^{n\times n}}
\le\frac{\|A_1-A_2\|_{\bbC^{n\times n}} 
+\|B_1-B_2\|_{\bbC^{n\times n}}}{1-2\|A_1\|_{\bbC^{n\times n}}}.
\end{equation}
Analogously, one proves that
\begin{equation}\lb{j10}
\|X_1-X_2\|_{\bbC^{n\times n}}
\le\frac{\|A_1-A_2\|_{\bbC^{n\times n}} 
+\|B_1-B_2\|_{\bbC^{n\times n}}}{1-2\|A_2\|_{\bbC^{n\times n}}}.
\end{equation}
Combining \eqref{j9} and \eqref{j10} yields \eqref{j8}.
\end{proof}

After these preparations we are able to prove the following 
fundamental uniqueness result.

\begin{theorem} \lb{t3.5}
Assume Hypothesis~\ref{h2.2}   
and let $x_0\in\bbR$. \\
$(i)$ In the case of a matrix-valued 
Schr\"odinger operator $H$, $g(z,x_0)$ and
$g'(z,x_0)$ for all $z\in\bbC_+$, uniquely 
determine the matrix-valued  $m\times m$ 
potential $Q(x)$ for a.e. $x\in\bbR$. \\
$(ii)$ In the case of a 
Dirac-type operator $D$ assume $x_0$ is a point of Lebesgue 
continuity of $B$ and define $g'(z,x_0)$ by the right-hand side of
\eqref{3.40}. Then $g(z,x_0)$, $g'(z,x_0)$ for all $z\in\bbC_+$,
and  the $m\times m$ coefficients $B_{1,1}(x_0)$, 
$B_{1,2}(x_0)$ uniquely determine the $2m\times 2m$ 
matrix $B(x)=\left(\begin{smallmatrix}B_{1,1}(x)& B_{1,2}(x)
\\B_{1,2}(x) & -B_{1,1}(x)\end{smallmatrix}\right)$ 
for a.e.~$x\in\bbR$.
\end{theorem}
\begin{proof}
Since $M(z,x_0)$ uniquely determines $T$ (i.e., $H$ or $D$) 
by Theorem~\ref{t2.5}, it suffices to show that 
$g(z,x_0)$ and
$g'(z,x_0)$ uniquely determine $M(z,x_0)$ in 
 case (i) and similarly, $g(z,x_0)$, $g'(z,x_0)$, and 
$B_{1,1}(x_0)$, $B_{1,2}(x_0)$ uniquely determine 
$M(z,x_0)$ in case (ii). We start with the
Schr\"odinger case (i). By \eqref{2.22}, \eqref{3.3}, and
\eqref{3.4}, 
\begin{equation}
\Im(g(z,x_0)) > 0 \text{ for all } z\in\bbC_+ \lb{3.48}
\end{equation}
and hence 
\begin{equation}
\spec(g(z,x_0))\cap\spec(-g(z,x_0))=\emptyset 
\text{ for all } z\in\bbC_+
\lb{3.49}
\end{equation}
by \eqref{3.45a}. Applying Lemma~\ref{l3.3} (with 
$A=-B=g(z,x_0)$ and $C=g'(z,x_0)\mp I_n$) to 
\eqref{3.39} (with $x=x_0$) then shows that $M_\pm(z,x_0)$ are
uniquely  determined by $g(z,x_0)$ and $g'(z,x_0)$. 

In the corresponding Dirac case (ii) one rewrites \eqref{3.40} 
(with $x=x_0$) in the form
\begin{align}
&g(z,x_0)^{-1}[I_m -z^{-1}B_{1,1}(x_0)]M_+(z,x_0)
+M_+(z,x_0)[I_m-z^{-1}B_{1,1}(x_0)]g(z,x_0)^{-1} 
 \no \\
&=z^{-1}[g(z,x_0)^{-1}g'(z,x_0)g(z,x_0)^{-1}-
B_{1,2}(x_0)g(z,x_0)^{-1}-g(z,x_0)^{-1}B_{1,2}(x_0) \no \\
& \quad -zg(z,x_0)^{-2}+g(z,x_0)^{-1}B_{1,1}(x_0)g(z,x_0)^{-1}].
\lb{3.50}
\end{align}
Thus, for $\Im(z)>0$ sufficiently large, we can again apply 
Lemma~\ref{l3.3} (with $A=g(z,x_0)^{-1}[I_m -
z^{-1}B_{1,1}(x_0)]$, $B=-[I_m-z^{-1}B_{1,1}(x_0)]g(z,x_0)^{-1}$, 
and $C=$ the r.-h.s. of \eqref{3.50}) since clearly
\begin{align}
&\spec (g(z,x_0)^{-1}[I_m - z^{-1}B_{1,1}(x_0)])\cap\spec
(-[I_m-z^{-1}B_{1,1}(x_0)]g(z,x_0)^{-1}) = \emptyset 
\lb{3.51} 
\end{align}
for $\Im(z)>0$ sufficiently large. 
Hence $M_+(z,x_0)$ is uniquely determined for $\Im(z)>0$
sufficiently large. Since $M_+(\cdot,x_0)$ is analytic in 
$\bbC_+$, one concludes that $M_+(z,x_0)$ is uniquely 
determined for all $z\in\bbC_+$. Since
$g(z,x_0)=[M_-(z,x_0)-M_+(z,x_0)]^{-1}$ is known by 
hypothesis, this also uniquely determines $M_-(z,x_0)$ 
for all $z\in\bbC_+$.

In either case an application of Theorem~\ref{t2.5} completes 
the proof.
\end{proof}

Since $g(z,x_0)$ and $g'(z,x_0)$ are analytic in $\bbC_+$, 
it suffices of course in Theorem~\ref{t3.5} to know them 
for an infinite set of $z\in\bbC_+$ with an accumulation 
point in $\bbC_+$.

In the special scalar case $m=1$, the uniqueness result 
stated in Theorem~\ref{t3.5} is well-known for
one-dimensional Schr\"odinger operators on $\bbR$ and can 
be inferred, for instance,  from the results in \cite{JM82}, 
\cite{Ko84}, \cite{Ko85}, \cite{Ko87} (it is also explicitly 
stated in \cite{GS96}). The present matrix extension 
with $m>1$, however, appears to be new. To the best of our
knowledge, the corresponding Dirac case has not previously 
been studied in the literature from this uniqueness point 
of view.

Next we turn to a local version of Theorem~\ref{t3.5}, 
our principal new result.

\begin{theorem} \lb{t3.5a}
Suppose $B_j$, $j=1,2$ satisfy Hypothesis~\ref{h2.2} and 
let $a>0$ and $x_0\in\bbR$. Introduce
$g_j(z,x_0)=  [M_{j,-}(z,x_0)-M_{j,+}(z,x_0)]^{-1}$, where
$M_{j,\pm}(z,x_0)$,
$j=1,2$  denote the $m\times m$ Weyl-Titchmarsh matrices associated 
with the corresponding Hamiltonian systems \eqref{2.4} on 
$(x_0,\pm\infty)$. \\
$(i)$ In the case of matrix-valued Schr\"odinger operators 
$H_j$, $j=1,2$ suppose that 
\begin{equation}
\|g_1(z,x_0)-g_2(z,x_0)\|_{\bbC^{m\times m}}+
\|g_1^\prime (z,x_0)-g_2^\prime (z,x_0)\|_{\bbC^{m\times m}}
\underset{|z|\to\infty}{=}O\big(e^{-2\Im(z^{1/2})a}\big) 
\lb{3.51a}
\end{equation}
along a ray with $0<\arg (z)<\pi$. Then
\begin{equation}
Q_1(x)=Q_2(x) \text{ for a.e. } x\in [x_0-a,x_0+a]. \lb{3.51b}
\end{equation}
$(ii)$ In the case of Dirac-type operators assume that $x_0$ 
is a point of Lebesgue continuity of $B_j$ and define
$g'_j(z,x_0)$, $j=1,2$, by the right-hand side of
\eqref{3.40}. Moreover, suppose that 
\begin{equation}
\|g_1(z,x_0)-g_2(z,x_0)\|_{\bbC^{m\times m}}+
\|g_1^\prime (z,x_0)-g_2^\prime (z,x_0)\|_{\bbC^{m\times m}}
\underset{|z|\to\infty}{=}O\big(e^{-2\Im(z)a}\big) 
\lb{3.51d}
\end{equation}
along a ray with $0<\arg (z)<\pi/2$ and along a ray with 
$\pi/2<\arg(z)<\pi$. In addition, assume that 
\begin{equation}
\big(B_1(x_0)\big)_{1,1}=\big(B_2(x_0)\big)_{1,1}, \quad 
\big(B_1(x_0)\big)_{1,2}=\big(B_2(x_0)\big)_{1,2}. \lb{3.51c}
\end{equation}
Then
\begin{equation}
B_1(x)=B_2(x) \text{ for a.e.~$x\in [x_0-a,x_0+a]$}. \lb{3.51e}
\end{equation}
\end{theorem}
\begin{proof}
In the Schr\"odinger context one observes from 
Lemmas~\ref{l3.3} and \ref{l3.4}, and \eqref{3.39} (for $x=x_0$)
that
\begin{align}
&M_{1,\pm}(z,x_0)-M_{2,\pm}(z,x_0) \no \\
&=(2\pi i)^{-1}\oint_{\Gamma_{1,2}} 
d\zeta\,\big[(g_1(z,x_0)-\zeta)^{-1}[g^\prime_1(z,x_0)\mp I_m]
(-g_1(z,x_0)-\zeta)^{-1} \no \\
&\hspace*{3cm} -(g_2(z,x_0)-\zeta)^{-1}[g^\prime_2 (z,x_0) \mp I_m]
(-g_2(z,x_0)-\zeta)^{-1}\big]. \lb{3.51f}
\end{align}
Here $\Gamma_{1,2}$ is a closed and counterclockwise oriented
Jordan contour which encircles $\spec(g_1(z,x_0))\cup\spec
(g_2(z,x_0))$ with winding number
$+1$ and $\spec (-g_1(z,x_0))\cup\spec (-g_2(z,x_0))$ with winding
number $0$, and whose existence is guaranteed by Lemma~\ref{l3.4}
(cf.~\eqref{3.49}). Using the asymptotic expansion 
(cf.~\cite{CG99})
\begin{equation}
g_j(z,x_0)\underset{|z|\to\infty}{=}(i/2)z^{-1/2}I_m+o(z^{-1}), 
\quad 0<\arg(z)<\pi, \; j=1,2, \lb{3.51g} 
\end{equation}
\eqref{3.51a} and a standard $\varepsilon/3$-argument yield
\begin{equation}
\|M_{1,\pm}(z,x_0)-M_{2,\pm}(z,x_0)\|_{\bbC^{m\times
m}}\underset{|z|\to\infty}{=}
O\big(e^{-2\Im(z^{1/2})a}\big) \lb{3.51h}
\end{equation}
along a ray with $\arg(z)=\pi-\varepsilon$. Hence one 
obtains \eqref{3.51b} by Theorem~\ref{t2.6}\,(i). 

In the case of Dirac-type operators one infers from 
\eqref{3.50} that 
\begin{align}
&g_j(z,x_0)^{-1}[I_m -z^{-1}\big(B_{j}(x_0)\big)_{1,1}]
M_{j,+}(z,x_0) \no \\
& \quad +M_{j,+}(z,x_0)[I_m-z^{-1}\big(B_{j}(x_0)\big)_{1,1}]
g_j(z,x_0)^{-1} \no \\
&=z^{-1}[g_j(z,x_0)^{-1}g_j^\prime(z,x_0)g_j(z,x_0)^{-1}-
\big(B_{j}(x_0)\big)_{1,2}g_j(z,x_0)^{-1} \no \\
& \quad -g_j(z,x_0)^{-1}\big(B_{j}(x_0)\big)_{1,2} 
-zg_j(z,x_0)^{-2} 
+g_j(z,x_0)^{-1}\big(B_{j}(x_0)\big)_{1,1}
g_j(z,x_0)^{-1}, \no \\
&\hspace*{10cm} j=1,2 \lb{3.51i}
\end{align}
and then proceeds in exactly the same manner as in 
\eqref{3.51f}--\eqref{3.51h} replacing \eqref{3.51g} by
(cf.~\cite{CG01})
\begin{equation}
g_j(z,x_0)\underset{|z|\to\infty}{=}(i/2)I_m+o(1), 
\quad 0<\arg(z)<\pi, \; j=1,2. \lb{3.52} 
\end{equation}
\end{proof}

For an alternative proof of \eqref{3.51h} one can multiply
\eqref{3.39} (for $x=x_0$) and \eqref{3.51i} by $(-i)$ and directly
apply  Lemma~\ref{l3.3}\,(ii) taking into account \eqref{3.51g} 
and \eqref{3.52}, respectively.

As mentioned in \cite{GS99}, the ray with 
$0<\arg(z)<\pi$ in Theorem~\ref{t3.5a}\,(i) is of no
particular significance and can be replaced by any
non-selfintersecting  Jordan arc that tends to infinity in the
sector $\varepsilon\leq\arg(z)\leq\pi-\varepsilon$ for some  
$0<\varepsilon <\pi/2$. The analogous comment applies to 
Theorem~\ref{t3.5a}\,(ii)

\vspace*{2mm}

As discussed in Marchenko's monograph
\cite[p.~193--195]{Ma86},  it is possible to reduce certain Dirac
operators $D$ with potentials $B(x)$ of the type
\eqref{2.13a} to $2m\times 2m$ matrix-valued Schr\"odinger 
operators by essentially  taking the square of the Dirac 
operator D. This permits one to relate the $M$-functions of the 
$2m\times 2m$ Schr\"odinger and Dirac-type operators and we next 
explore this connection in some detail. In order to indicate that
$B$ is temporarily of the special form \eqref{2.13a}, we denote it
by
$\wti B$ in the  following and assume in addition to
Hypothesis~\ref{h2.2} that
\begin{equation}
\wti B=\begin{pmatrix}
\wti B_{1,1}& \wti B_{1,2} 
\\ \wti B_{1,2} & -\wti B_{1,1}
\end{pmatrix}\in C^1(\bbR)^{2m\times 2m} \lb{3.53} 
\end{equation}
(this can be improved but is of little relevance
in our  context). Denoting the resulting Dirac operator by $\wti D$,
\begin{equation}
\wti D=J\f{d}{dx}-\wti B, \lb{3.54}
\end{equation}
one computes, using $J^2=-I_{2m}$ and $J\wti B(x)+\wti B(x)J=0$, 
\begin{equation}
{\wti D}^2 = -I_{2m}\f{d^2}{dx^2} + {\wti B}^2-J\wti B' 
\lb{3.55}
\end{equation}
and hence obtains the $2m\times 2m$ matrix-valued Schr\"odinger
operator 
\begin{equation}
\wti H = -I_{2m}\f{d^2}{dx^2} + \wti Q, \lb{3.56}
\end{equation}
where
\begin{align}
\wti Q (x)&= \wti B(x)^2-J\wti B'(x)  \lb{3.57} \\
&= \left(\begin{smallmatrix}
\wti B_{1,1}(x)^2+\wti B_{1,2}(x)^2+\wti B_{1,2}'(x)&
\wti B_{1,1}(x)\wti B_{1,2}(x)-\wti B_{1,2}(x)\wti B_{1,1}(x)-\wti
B_{1,1}'(x) 
\\ \wti B_{1,2}(x)\wti B_{1,1}(x)-\wti B_{1,1}(x)\wti
B_{1,2}(x)-\wti B_{1,1}'(x) & \wti B_{1,1}(x)^2+\wti
B_{1,2}(x)^2-\wti B_{1,2}'(x)
\end{smallmatrix}\right). \no 
\end{align}
Since
\begin{equation}
(\wti H-z^2)=(\wti D-z)(\wti D+z), \lb{3.58}
\end{equation}
linearly independent solutions $\Phi_+ (z^2,x,x_0),
\Phi_- (z^2,x,x_0)$ of 
$(\wti H -z^2) \Phi (z^2,x,x_0)=0$ are of the type
\begin{equation}
\Phi_\pm(z^2,x,x_0)=\begin{pmatrix}
\psi_{\pm,1}(z,x,x_0)& \psi_{\pm,1}(-z,x,x_0)  \\ 
\psi_{\pm,2}(z,x,x_0) & \psi_{\pm,2}(-z,x,x_0)
\end{pmatrix}, 
\quad z^2\in\bbC_+, \,\,x,x_0\in\bbR, \lb{3.59}
\end{equation}
where, in accordance with \eqref{2.31}, the $m$ columns of
$\Psi_\pm(z,x,x_0)$
\begin{align}
\Psi_\pm(z,x,x_0)&=\begin{pmatrix}\psi_{\pm,1}(z,x,x_0)\\
\psi_{\pm,2}(z,x,x_0)  \end{pmatrix}  \lb{3.60} \\
&=\begin{pmatrix}\theta_1(z,x,x_0) & \phi_1(z,x,x_0)\\
\theta_2(z,x,x_0)& \phi_2(z,x,x_0)\end{pmatrix} 
\begin{pmatrix} I_m \\
M_\pm^{\wti D} (z,x_0) \end{pmatrix},  
\quad z\in\bbC\backslash\bbR, \,\, x,x_0\in\bbR \no
\end{align}
form a basis for the nullspace of $(\wti D -z)$ Here, in 
obvious notation, $M_\pm^{\wti D} (z,x_0)$ abbreviates the 
corresponding $m\times m$ half-line Weyl-Titchmarsh matrices
associated  with $\wti D$, given by
\begin{align}
M_\pm^{\wti D} (z,x_0)&=\psi_{\pm,2}(z,x_0,x_0)
\psi_{\pm,1}(z,x_0,x_0)^{-1}=\psi_{\pm,2}(z,x_0,x_0),  \no \\
M_\pm^{\wti D} (-z,x_0)&=M_\pm^{\wti D} (-\overline z,x_0)^*, 
\quad z\in\bbC_+, \lb{3.61}
\end{align}
using the normalization $\psi_{\pm,j}(z,x_0,x_0)=I_m$, 
$j=1,2$. The corresponding $2m\times 2m$ half-line 
Weyl-Titchmarsh matrices $M_\pm^{\wti H} (z^2,x_0)$, associated 
with $\wti H$, are then defined by
\begin{equation}
M_\pm^{\wti H} (z^2,x_0)=\Phi_\pm (z^2,x_0,x_0)
\Phi_\pm (z^2,x_0,x_0)^{-1}, \quad
z^2\in\bbC_+. \lb{3.62}
\end{equation}
Next we explicitly compute $M_\pm^{\wti H} (z^2,x_0)$ in terms 
of $M_\pm^{\wti D} (z,x_0)$, $M_\pm^{\wti D} (-z,x_0)$.
\begin{lemma} \lb{l3.6} Assume Hypothesis~\ref{h2.2},  
\eqref{3.53}, and $\arg (z)\in (0,\pi/2)$. 
Then
\begin{align}
M_\pm^{\wti H} (z^2,x_0)&= \big(M_\pm^{\wti H}
(z,x_0)_{j,k}\big)_{1\leq j,k\leq 2}, \lb{3.63} \\
M_\pm^{\wti H} (z^2,x_0)_{1,1}&=
\wti B_{1,2}(x_0)\no \\
& \quad +2zM_\pm^{\wti D} (z,x_0)[M_\pm^{\wti D}
(-z,x_0)-M_\pm^{\wti D} (z,x_0)]^{-1}M_\pm^{\wti D} (-z,x_0), \no \\
M_\pm^{\wti H}(z^2,x_0)_{1,2}&=-\wti B_{1,1}(x_0) \no \\
& \quad -z[M_\pm^{\wti D} 
(-z,x_0)+M_\pm^{\wti D} (z,x_0)][M_\pm^{\wti D} 
(-z,x_0)-M_\pm^{\wti D} (z,x_0)]^{-1}, \no \\
M_\pm^{\wti H}(z^2,x_0)_{2,1}&=-\wti B_{1,1}(x_0) \no \\
& \quad -z[M_\pm^{\wti D} 
(-z,x_0)+M_\pm^{\wti D} (z,x_0)]^{-1}[M_\pm^{\wti D} 
(-z,x_0)-M_\pm^{\wti D} (z,x_0)], \no \\
M_\pm^{\wti H}(z^2,x_0)_{2,2}&=-\wti B_{1,2}(x_0)+2z[M_\pm^{\wti D} 
(-z,x_0)-M_\pm^{\wti D} (z,x_0)]^{-1}. \no 
\end{align}
\end{lemma}
\begin{proof}
Combining \eqref{3.62}, 
\begin{align}
&\Phi_\pm(z^2,x_0,x_0)=\begin{pmatrix}
I_m & I_m \\ M_\pm^{\wti D} (z,x_0) & M_\pm^{\wti D} (-z,x_0)
\end{pmatrix}, \lb{3.64} \\
&\Phi_\pm(z^2,x_0,x_0)^{-1} \lb{3.65} \\
&=\begin{pmatrix}
[M_\pm^{\wti D}  
(-z,x_0)-M_\pm^{\wti D} (z,x_0)] M_\pm^{\wti D} (-z,x_0) & 
-[M_\pm^{\wti D} (-z,x_0)-M_\pm^{\wti D} (z,x_0)]^{-1} \\ 
-[M_\pm^{\wti D} (-z,x_0)-M_\pm^{\wti D} (z,x_0)] 
M_\pm^{\wti D} (z,x_0)&  [M_\pm^{\wti D} 
(-z,x_0)-M_\pm^{\wti D} (z,x_0)]^{-1}
\end{pmatrix}, \no \\
&J\Psi_\pm^{\prime}(z,x_0,x_0)=\begin{pmatrix}
\wti B_{1,1}(x_0)+zI_m + \wti B_{1,2}(x_0)M_\pm^{\wti D}(z,x_0) \\ 
\wti B_{1,2}(x_0) -\wti B_{1,1}(x_0)M_\pm^{\wti D}(z,x_0)
+zM_\pm^{\wti D}(z,x_0)
\end{pmatrix}, \lb{3.66}
\end{align}
yields \eqref{3.63} after straightforward matrix computations.
\end{proof}

\section{Matrix-Valued Jacobi Operators} \lb{s4}

In this section we prove analogous uniqueness theorems for 
matrix-valued Jacobi (i.e., second-order finite difference)
operators. We could consider analogous  Dirac-type difference
operators and hence introduce discrete Hamiltonian systems as 
in Section~\ref{s2}, but for simplicity we will confine 
ourselves to the Jacobi case.

As the basic hypothesis of this section we adopt the following.

\begin{hypothesis}\lb{h4.1}
Let $m\in\bbN$ and consider the sequences of 
self-adjoint $m\times m$ matrices
\begin{align}
&A=\{A(k)\}_{k\in\bbZ}, \quad 
A(k)=A(k)^*, \,\, A(k)>0, \,\, k\in\bbZ, \no \\
&B=\{B(k)\}_{k\in\bbZ}, \quad B(k)=B(k)^*, \,\, k\in\bbZ. \lb{4.1}
\end{align}
Moreover, we assume that $A(k)$ and $B(k)$ are uniformly bounded
with  respect to $k$, that is, there exists a $C>0$,
such that 
\begin{equation}
\|A(k)\|_{\bbC^{m\times m}}+\|B(k)\|_{\bbC^{m\times m}}\leq C,
\quad k\in\bbZ. \lb{4.1a}
\end{equation}
\end{hypothesis}
Next, denote by $S^\pm$ the shift operators in 
$\ell^\infty (\bbZ)^m$ and $\ell^\infty (\bbZ)^{m\times m}$, 
that is, 
\begin{equation}
(S^\pm g)(k)=g^\pm (k)=g(k\pm 1) \text{ for } g\in \ell^\infty
(\bbZ)^m \text{ or }  g\in \ell^\infty (\bbZ)^{m\times m}, \,\,
k\in\bbZ.  \lb{4.7}
\end{equation}
Given Hypothesis~\ref{h4.1}, the matrix-valued Jacobi 
operator $H$ in $\ell^2(\bbZ)^m$ is then defined by
\begin{equation} 
H = A S^+ + A^- S^-  + B, \quad \dom(H)=\ell^2(\bbZ)^m.  
\label{4.8b} 
\end{equation}
Because of hypothesis \eqref{4.1a}, $H$ is a bounded symmetric
operator and hence self-adjoint. In particular, the difference 
expression $A S^+ + A^- S^-  + B$ induced by \eqref{4.8b} is in the
limit point case at $\pm\infty$. We chose to adopt \eqref{4.1a} for
simplicity only.  Our formalism extends to unbounded Jacobi
operators and to cases  where the difference expression associated
with \eqref{4.8b} is in  the limit circle case at $+\infty$ and/or
$-\infty$.  We just note in passing that without assuming
\eqref{4.1a}, the difference expression $A S^+ + A^- S^-  + B$ is in
the limit point case at
$\pm\infty$ if $\sum_{k=k_0}^{\pm\infty}\|A(k)\|_{\bbC^{m\times
m}}^{-1} =\infty$ (see, e.g., \cite[Theorem~VII.2.9]{Be68}). 

The Green's matrix associated
with $H$ will be denoted by $G(z,k,\ell)$ in the following,
\begin{equation}
G(z,k,\ell)=(H-z)^{-1}(k,\ell), \quad z\in\bbC\backslash\bbR, \,\, 
k,\ell\in\bbZ. \lb{4.8a}
\end{equation}

Next, fix a site $k_0 \in \bbZ$ and in analogy to 
\eqref{2.15}--\eqref{2.17} define $m
\times m$  matrix-valued solutions $\phi(z,k,k_0)$ and
$\theta(z,k,k_0)$ of the equation
\begin{equation}
A(k)\psi(z,k+1)+A(k-1)\psi(z,k-1)+(B(k)-zI_m)\psi(z,k)=0, 
\quad z\in\bbC, \, k\in\bbZ, \lb{4.2}
\end{equation}
satisfying the initial conditions
\begin{equation}\label{4.3}
\theta(z,k_0,k_0)=\phi(z,k_0+1,k_0)=I_m, \quad 
\phi(z,k_0,k_0)=\theta(z,k_0+1,k_0)=0.
\end{equation}
As in \eqref{2.31}, \eqref{2.32} one then introduces 
$m\times m$ matrix-valued Weyl solutions $\psi_{\pm}(z,k,k_0)$
associated with $H$ defined by 
\begin{equation}
\psi_{\pm}(z,k,k_0)= \theta (z,k,k_0) - \phi (z,k,k_0)A(k_0)^{-1}
M_\pm(z,k_0), \quad z\in\bbC\backslash\bbR, \,\, k\in\bbZ, \lb{4.4}
\end{equation}
with the properties 
\begin{align}
&A(k)\psi_{\pm}(z,k+1,k_0)+A(k-1)\psi_\pm(z,k-1,k_0)
+(B(k)-z)\psi_\pm(z,k,k_0)=0, \no \\
& \hspace*{8.4cm} z\in\bbC\backslash\bbR, \,\, k\in\bbZ, 
\lb{4.4a} \\ 
&\psi_\pm (z,\cdot,k_0)\in
\ell^2((k_0,\pm\infty)\cap
\bbZ)^{m\times m}, \quad z\in\bbC\backslash\bbR, \lb{4.5} 
\end{align}
where $M_\pm(z,k_0)$ denote the Weyl-Titchmarsh matrices
associated with $H$. Since by assumption $A S^+ + A^- S^- + B$ 
is in the limit point case at $\pm\infty$, $M_\pm(z,k_0)$ in
\eqref{4.4} are uniquely determined by the requirement \eqref{4.5}.
We also note that by a standard argument, 
\begin{equation}
\det(\phi(z,k,k_0))\neq 0 \text{ for all }
k\in\bbZ\backslash\{k_0\} \text{ and } z\in\bbC\backslash\bbR, 
\lb{4.4b}
\end{equation}
since otherwise one could construct a Dirichlet-type
eigenvalue $z\in\bbC\backslash\bbR$ for $H$ restricted to the finite
segment $k_0+1,\dots,k-1$ for
$k\geq k_0+1$ of $\bbZ$ (and similarly for $k\leq k_0-1$).
Thus, introducing 
\begin{equation}
M_N(z,k_0)=-\phi(z,N,k_0)^{-1}\theta(z,N,k_0), \quad
z\in\bbC\backslash\bbR, \,\, N\in\bbZ\backslash\{k_0\}, \lb{4.5a}
\end{equation}
one can then compute $M_\pm(z,k_0)$ by the limiting relation
\begin{equation}
M_\pm(z,k_0)=\lim_{N\to\pm\infty} M_N(z,k_0), \quad
z\in\bbC\backslash\bbR, \lb{4.5b}
\end{equation}
the limit being unique since $A S^+ + A^- S^-  + B$ is in
the limit point case at $\pm\infty$. Alternatively, \eqref{4.4}
yields
\begin{equation}
M_\pm(z,k_0)=-A(k_0)\psi_\pm(z,k_0+1,k_0), \quad
z\in\bbC\backslash\bbR. \lb{4.5c}
\end{equation}
More generally, recalling
\begin{equation}
\det(\psi_\pm(z,k,k_0))\neq 0 \text{ for all }
k\in\bbZ \text{ and } z\in\bbC\backslash\bbR, 
\lb{4.5d}
\end{equation}
by an argument analogous to that following \eqref{4.4b}, we now 
introduce
\begin{equation}
M_\pm(z,k)=-A(k)\psi_\pm(z,k+1,k_0)\psi_\pm(z,k,k_0)^{-1}, \quad
z\in\bbC\backslash\bbR, \,\, k\in\bbZ. \lb{4.5e}
\end{equation}
One easily verifies that $M_\pm(z,k)$ represent the
Weyl-Titchmarsh $M$-matrices associated with the reference point
$k\in\bbZ$. Moreover, one obtains the Riccati-type equation for
$M_\pm(z,k)$, 
\begin{equation}
M_\pm(z,k)+A(k-1)M_\pm(z,k-1)^{-1}A(k-1)-(B(k)-z)=0, \quad 
z\in\bbC\backslash\bbR, \,\, k\in\bbZ \lb{4.5f}
\end{equation}
as a result of \eqref{4.4a}.

Next, we define the self-adjoint half-line Jacobi operators 
$H_{\pm,k_0}$ in $\ell^2((k_0,\pm\infty)\cap\bbZ)^m$ by 
\begin{equation}\lb{4.9}
H_{\pm, k_0}=P_{\pm, k_0}HP_{\pm,
k_0}\big|_{\ell^2((k_0,\pm\infty)\cap\bbZ)^m}, \quad k_0\in \bbZ,
\end{equation}
where $P_{\pm, k_0}$ are the orthogonal
projections onto the subspaces
$\ell^2([k_0,\pm \infty)\cap \bbZ)^m$. By inspection, $H_{\pm,k_0}$
satisfy Dirichlet boundary conditions at $k_0\mp 1$, 
\begin{align} 
(H_{+, k_0}f)(k_0) &= A(k_0)f^+(k_0)+B(k_0)f(k_0), \no \\
(H_{-, k_0}f)(k_0) &= A^-(k_0)f^-(k_0)+B(k_0)f(k_0), \no \\
(H_{\pm, k_0}f)(k) &= A(k)f^+(k)+A^-(k)f^-(k)+B(k)f(k), \quad 
k\gtreqless k_0\pm 1, \label{4.8}  \\
f\in \dom(H_{\pm, k_0})&=\ell^2((k_0,\pm\infty)\cap\bbZ)^m \no
\end{align}
(i.e., formally, $f(k_0\mp 1)=0$). We also introduce Weyl
$m$-matrices $m_\pm(z,k_0)$ associated with $H_{\pm, k_0}$ by
\begin{align}
m_\pm(z,k_0)&=Q_{k_0}(H_{\pm, k_0}-z)^{-1}Q_{ k_0}, \lb{4.10} \\
&=G_{\pm,k_0}(z,k_0,k_0), \quad z\in\bbC\backslash\bbR,
\,\, k_0\in \bbZ. \lb{4.10a}
\end{align}
Here $Q_{ k_0}$ are orthogonal projections
onto the $m$-dimensional subspaces $\ell^2(\{k_0 \})^m$,
$k_0\in\bbZ$  and 
\begin{equation}
G_{\pm,k_0}(z,k,\ell)=(H_{\pm,k_0}-z)^{-1}(k,\ell), \quad
z\in\bbC\backslash\bbR, \,\,  k,\ell\in\bbZ \lb{4.10b}
\end{equation}
represents the Green's matrix of $H_{\pm,k_0}$.

In order to find the connection between $m_\pm(z,k_0)$ and
$M_\pm(z,k_0)$ we briefly discuss the Green's matrices
$G_{\pm,k_0}(z,k,\ell)$ and $G(z,k,\ell)$ associated
with $H_{\pm, k_0}$ and $H$ next. 

First we recall the definition of the Wronskian $W(f,g)(k)$ of two
sequences of matrices $f(\cdot), g(\cdot)\in \bbC^{m\times m} $
given by 
\begin{equation}\lb{4.6}
W(f,g)(k)=f(k)A(k)g(k+1)-f(k+1)A(k)g(k), \quad k\in\bbZ.
\end{equation}
We note that for any two matrix-valued solutions 
$\varphi(z,\cdot)$ and  $\psi(z,\cdot)$ of \eqref{4.2}
the Wronskian $W(\varphi(\overline{ z},\cdot)^*,\psi( z,\cdot))(k)$
is independent of $k\in\bbZ$.

In complete analogy to the scalar Jacobi case (i.e., $m=1$) one
verifies,
\begin{align}
&G_{+, k_0}(z,k,\ell)=\begin{cases} -\psi_+(z,k,k_0-1) 
A(k_0-1)^{-1}\phi(\ol z, \ell,k_0-1)^*, & \ell\leq k, \\
-\phi(z,k,k_0-1) 
A(k_0-1)^{-1}\psi_+(\ol z, \ell,k_0-1)^*, & \ell\geq k, 
\end{cases} \lb{4.10c} \\
&G_{-, k_0}(z,k,\ell)=\begin{cases} \phi(z,k,k_0+1) 
A(k_0+1)^{-1}\psi_-(\ol z, \ell,k_0+1)^*, & \ell\leq k, \\
\psi_-(z,k,k_0+1) 
A(k_0+1)^{-1}\phi(\ol z, \ell,k_0+1)^*, & \ell\geq k, 
\end{cases} \lb{4.10c'} \\
& \hspace*{7.7cm} z\in\bbC\backslash\bbR, \,\,  k,\ell\in\bbZ. \no
\end{align}
Similarly, using the fact that 
\begin{align}
&\psi_+(z,k,k_0)[M_-(z,k_0)
-M_+(z,k_0)]^{-1}\psi_- (\overline {z},k,k_0)^* \no \\
&= \psi_-(z,k,k_0)[M_-(z,k_0)
-M_+(z,k_0)]^{-1}\psi_+ (z,k,k_0) \no \\
&=[M_-(z,k)-M_+(z,k)]^{-1}, \quad k\in \bbZ, \lb{4.16}
\end{align}
and 
\begin{align}
&A(k)\psi_+(z,k+1,k_0)[M_-(z,k_0)
-M_+(z,k_0)]^{-1}\psi_- (\overline z,k,k_0)^*  \lb{4.17} \\
&- A(k)\psi_-(z,k+1,k_0)[M_-(z,k_0)-M_+(z,k_0)]^{-1}
\psi_+ (\overline z,k,k_0)^*=I_m, \quad k\in \bbZ, \no
\end{align}
one verifies 
\begin{align}
&G(z,k,\ell) = \begin{cases} \psi_+(z,k,k_0)[M_-(z,k_0)
-M_+(z,k_0)]^{-1}\psi_- (\overline z,\ell,k_0)^*,  & \ell \leq k,
 \\
\psi_-(z,k,k_0)[M_-(z,k_0)
-M_+(z,k_0)]^{-1}\psi_+ (\overline z,\ell,k_0)^*,  & \ell \geq k,
\end{cases} \label{4.15} \\
& \hspace*{8.3cm} z\in\bbC\backslash\bbR, \,\,  k,\ell\in\bbZ. \no
\end{align}

Using \eqref{4.3}, \eqref{4.4a}, \eqref{4.5f}, and \eqref{4.10a}, 
one infers that the Weyl-Titchmarsh matrices $M_\pm(z,k_0)$
introduced by \eqref{4.5b} (resp.~\eqref{4.5c}) and the Weyl
$m$-functions $m_\pm(z,k_0)$ defined in \eqref{4.10} are 
related by
\begin{equation} \lb{4.11}
M_+(z,k_0)=-m_+(z, k_0)^{-1}-zI_m+B(k_0)
\end{equation}
and
\begin{equation}\lb{4.12}
M_-(z,k_0)=m_-(z, k_0)^{-1}.
\end{equation}

As in the continuous case, 
$\pm M_\pm(z,k_0)$ are Herglotz  matrices satisfying 
\begin{equation}
\pm \Im (M_\pm (z,k_0))>0, \,\, z\in\bbC_+, \quad M_\pm (\overline
z,k_0)  = M_\pm (z,k_0)^*. \lb{4.13}
\end{equation}

Next, we introduce the $2m\times 2m$ Weyl-Titchmarsh matrix 
$M(z,k_0)$ associated with the Jacobi operator $H$ in 
$\ell^2(\bbZ)^m$ as in \eqref{2.34} (with
$x_0\in\bbR$ replaced by $k_0\in\bbZ$, etc.) by
\begin{equation}
M(z,k_0)=\big(M_{j,j^\prime}(z,k_0)\big)_{j,j^\prime=1,2}\, , 
\quad z\in\bbC\backslash\bbR, \,\, k_0\in\bbZ, \lb{4.21}
\end{equation}
where
\begin{align}
M_{1,1}(z,k_0)&=[M_{-}(z,k_0)-M_{+}(z,k_0)]^{-1}, \lb{4.22} \\
M_{1,2}(z,k_0)&=2^{-1}[M_{-}(z,k_0)-M_{+}(z,k_0)]^{-1}
[M_{-}(z,k_0)+M_{+}(z,k_0)], \lb{4.23} \\
M_{2,1}(z,k_0)&=2^{-1}[M_{-}(z,k_0)+M_{+}(z,k_0)]
[M_{-}(z,k_0)-M_{+}(z,k_0)]^{-1}, \lb{4.23a} \\
M_{2,2}(z,k_0)&=M_{\pm}(z,k_0)
[M_{-}(z,k_0)-M_{+}(z,k_0)]^{-1}M_{\mp}(z,k_0). \lb{4.23b} 
\end{align}

\begin{remark}
We note that the Weyl-Titchmarsh matrix 
$M(z,k_0)$ is related to the Green's matrix associated with the 
Jacobi operator $H$ by
\begin{equation}
M(z,k_0)=\begin{pmatrix} I_m &0  \\
0 &-A(k_0) \end{pmatrix} \calM(z,k_0)
\begin{pmatrix}I_m &0  \\
0 &-A(k_0) \end{pmatrix}
+\frac{1}{2}\begin{pmatrix} 0 &I_m  \\
I_m &0 \end{pmatrix},
\end{equation}
where
\begin{equation}
\calM(z,k_0)=
\begin{pmatrix}G(z, k_0,k_0) &G(z,k_0,k_0+1)  \\
 G(z,k_0+1, k_0) & G(z, k_0+1, k_0+1) \end{pmatrix}.
\end{equation}
\end{remark}

With $M(z,k_0)$ defined in \eqref{4.21}, the analog of
Theorem~\ref{t2.5} holds. 

\begin{theorem} \lb{t4.2}
Assume Hypothesis~\ref{h4.1} and let $k_0\in\bbZ$. Then
the $2m\times 2m$ Weyl-Titchmarsh matrix 
$M(z,k_0)$ for all $z\in\bbC_+$ uniquely determines the 
Jacobi operator $H$ and hence $A=\{A(k)\}_{k\in\bbZ}$ and 
$B=\{B(k)\}_{k\in\bbZ}$.
\end{theorem}

Perhaps the simplest way to prove Theorem~\ref{t4.2} is to reduce
it  to knowledge of $M_\pm(z,k_0)$, and hence by
\eqref{4.11} and \eqref{4.12} to that of $m_\pm(z,k_0)$ 
 for all $z\in\bbC_+$. Using the standard construction of orthogonal matrix-valued
polynomials with respect to the normalized measure in
the Herglotz representation of $m_\pm(z,k_0)$,
\begin{equation}
 m_\pm(z,k_0)=\int_\bbR d\Omega_\pm(\lambda,k_0) \,
(\lambda-z)^{-1}, \quad z\in\bbC_+, \quad
\int_\bbR d\Omega_\pm(\lambda,k_0) =I_m, \lb{4.24}
\end{equation}
allows one
to reconstruct $A(k)$, $B(k)$, $k\in [k_0,\pm\infty)\cap\bbZ$
from the measures $d\Omega_\pm(\lambda,k_0)$ (cf., e.g., 
\cite[Section~VII.2.8]{Be68}).
More precisely,
\begin{align}
&A(k_0\pm k)=\int_\bbR\lambda P_{\pm,k}(\lambda, k_0)
 d\Omega_\pm(\lambda,k_0) \,
P_{\pm,k+1}(\lambda, k_0)^*,\no \\
&
 B(k_0\pm k)=\int_\bbR\lambda P_{\pm,k}(\lambda,k_0)
 d\Omega_\pm(\lambda,k_0)
\, P_{\pm,k}(\lambda,k_0)^*, \quad k\in\bbN_0, \lb{4.25}
\end{align} 
where $\{P_{\pm,k}(\lambda,k_0)\}_{k\in\bbN_0}$ is an orthonormal 
system of matrix-valued polynomials with respect to the spectral 
measure $d\Omega_\pm(\lambda,k_0)$, with $P_{\pm,k}(z,k_0)$ 
of degree $k$ in $z$, $P_{\pm, 0}(z,k_0)=I_m$. One verifies, 
\begin{align}
P_{+,k}(z,k_0)&=\phi(z,k_0+k,k_0-1), \lb{4.25a} \\
P_{-,k}(z,k_0)&=\theta(z,k_0-k,k_0), \quad k\in\{-1\}\cup\bbN_0,
\,\, k_0\in\bbZ, \,\, z\in\bbC, \lb{4.25b}
\end{align}
with $\phi(z,k,k_0)$ and $\theta(z,k,k_0)$ defined in 
\eqref{4.2}, \eqref{4.3}. 
\bigskip

In the following we denote by $\lfloor x\rfloor\in\bbZ$  the
largest integer  less or equal to $x\in\bbR$.

The local analog of Theorem~\ref{t4.2} then reads as follows. 

\begin{theorem} \lb{t4.2a}
Suppose $0< \varepsilon<\pi$.~Let $H_j$, $j=1,2$ be two Jacobi
operators in 
$\ell^2 (\bbZ)^m$, with coefficients $\{A_{j,k}\}_{k\in\bbZ}$,
$\{B_{j,k}\}_{k\in\bbZ}$ satisfying 
Hypothesis~\ref{h4.1}, and denote by $M_{j,\pm}(z,k_0)$
 the $m\times m$ Weyl-Titchmarsh 
matrices associated with $H_j$, $j=1,2$. If there 
exists an $N\in\bbN$, $N\geq 2$, such that 
\begin{equation}
\|M_{1,-}(z,k_0)-M_{2,-}(z,k_0)
\|_{\bbC^{m\times m}}\underset{|z|\to\infty}
{=}O(|z|^{-N}) \lb{4.26}
\end{equation} 
along a ray with $\arg(z)=\pi -\varepsilon$, then 
\begin{align}
A_{1}(k_0- k)&=A_{2}(k_0- k), \quad k=0, \, ...\, ,
\lfloor N/2 \rfloor-1,
\no \\
 \quad B_{1}(k_0- k)&=B_{2}(k_0- k), \quad k=0, \, ...\, ,
 \lfloor (N-1)/2 \rfloor. \lb{4.27}
\end{align}
If there exists an $N\in\bbN$, $N\geq 2$, such that 
\begin{equation}
\|M_{1,+}(z,k_0)-M_{2,+}(z,k_0)
\|_{\bbC^{m\times m}}\underset{|z|\to\infty}
{=}O(|z|^{-N}) \lb{4.27c}
\end{equation} 
along a ray with $\arg(z)=\pi -\varepsilon$ and
\begin{equation}
B_1(k_0)=B_2(k_0), \lb{4.27a}
\end{equation}
then 
\begin{align}
A_{1}(k_0+ k)&=A_{2}(k_0+ k), \quad k=0, \, ...\, ,
\lfloor N/2 \rfloor-1,
\no \\
 \quad B_{1}(k_0+ k)&=B_{2}(k_0+ k), \quad k=0, \, ...\, ,
 \lfloor (N-1)/2 \rfloor. \lb{4.27b}
\end{align}
\end{theorem}
\begin{proof} Using \eqref{4.11}, \eqref{4.12}, 
and the asymptotic representations
\begin{equation}\lb{4.28}
m_{j,\pm}(z,k_0)\underset{|z|\to\infty}{=}
-z^{-1}+O(|z|^{-2}), \quad j=1,2,
\end{equation}
 one infers
\begin{equation}\lb{4.29}
m_{1,\pm}(z,k_0)-m_{2,\pm}(z,k_0)
\underset{|z|\to\infty}{=}
O\big (|z|^{-2}(M_{1,\pm}(z,k_0)-M_{2,\pm}(z,k_0))\big ).
\end{equation}
Here assumption \eqref{4.27a} is explicitly needed in the
case of the $+$-sign in \eqref{4.29}. By hypothesis
\eqref{4.26},
\eqref{4.29} implies
\begin{equation}\lb{4.30}
m_{1,\pm}(z,k_0)-m_{2,\pm}(z,k_0)
\underset{|z|\to\infty}{=} O\big (|z|^{-(N+2)}),
\end{equation}
and hence
\begin{equation}
\int_\bbR  d\Omega_{1,\pm}(\lambda,k_0) \lambda^n =
\int_\bbR  d\Omega_{2,\pm}(\lambda,k_0) \lambda^n 
, \quad n=0,1, \, ...\, , N. \lb{4.32} 
\end{equation}

Introducing the  orthogonal polynomials
 $\{P_{j,\pm;k}(\lambda,k_0)\}_{k\in\bbN_0}$
 associated with the normalized spectral measures
  $d\Omega_{j,\pm}(\lambda,k_0)$, $j=1,2$,  
\eqref{4.32} implies that
\begin{align}
&\int_\bbR\lambda P_{1,\pm;k}(\lambda, k_0)
 d\Omega_{1,\pm}(\lambda,k_0) \,
P_{1,\pm;k+1}(\lambda, k_0)^*
\no \\
&=\int_\bbR\lambda P_{2,\pm;k}(\lambda, k_0)
 d\Omega_{2,\pm}(\lambda,k_0) \,
P_{2,\pm;k+1}(\lambda, k_0)^*, 
\quad k=0,  ..., \,\lfloor N/2\rfloor -1, 
\lb{4.33}
\end{align}
and
\begin{align}
&
\int_\bbR\lambda P_{1,\pm;k}(\lambda,k_0)
 d\Omega_{1,\pm}(\lambda,k_0)
\, P_{1,\pm;k}(\lambda,k_0)^* 
\no \\
&=
\int_\bbR\lambda P_{2,\pm;k}(\lambda,k_0)
 d\Omega_{2,\pm}(\lambda,k_0)
\, P_{2,\pm;k}(\lambda,k_0)^*,
\quad k=0,  ..., \,\lfloor (N-1)/2\rfloor,
 \lb{4.34}
\end{align} 
which  proves \eqref{4.27} using \eqref{4.25}.
\end{proof}

Of course, one can permit $\varepsilon=0$ and $\varepsilon=\pi$ in
Theorem~\ref{t4.2a} (and in Theorem~\ref{t4.6} below) as long as
$|z|>\max_{j=1,2} \|H_j\|$. 

\begin{remark} \lb{r4.4a} 
The apparent asymmetry with respect to the additional condition
\eqref{4.27a}  in connection with $M_{j,+}(z,k_0)$ is due to the
fact that for  fixed $k_0\in\bbZ$, the map
\begin{equation}
(m_+(z,k_0),B(k_0)) \mapsto M_+(z,k_0), \quad z\in\bbC 
\lb{4.34a}
\end{equation}
given by \eqref{4.11} is not injective and hence knowledge of 
$M_+(z,k_0)$ does not determine $m_+(z,k_0)$ uniquely. This
follows from the explicit cancellation of $B(k_0)$ in 
\eqref{4.11} as $|z|\to\infty$ 
in sharp contrast to the asymptotic behavior of $M_-(z,k_0)$ 
in \eqref{4.12}, as described explicitly in \eqref{4.38} and 
\eqref{4.39} below.
\end{remark}

Given these preliminaries and introducing the diagonal Green's
matrix by 
\begin{equation}\lb{4.35}
g(z,k)=G(z,k,k), \quad z\in\bbC\backslash\bbR, \,\, k\in \bbZ,
\end{equation}
we can now formulate the analog of the uniqueness result,
Theorem~\ref{t3.5}, for Jacobi operators.

\begin{theorem}\label{t4.4} 
Assume Hypothesis~\ref{h4.1} and let $k_0\in\bbZ$. Then any
of the following three sets of data\\ 
$(i)$ $g(z,k_0)$ and $G(z,k_0,k_0+1)$ for all $z\in\bbC_+$; \\ 
$(ii)$ $g(z,k_0)$ and $[G(z,k_0,k_0+1)+G(z,k_0+1,k_0)]$ 
for all $z\in\bbC_+$; \\
$(iii)$ $g(z,k_0)$, $g(z,k_0+1)$ for all $z\in\bbC_+$ 
and $A(k_0)$; \\
uniquely determines the matrix-valued Jacobi operator $H$ 
and hence $A=\{A(k)\}_{k\in\bbZ}$ and $B=\{B(k)\}_{k\in\bbZ}$.\\ 
\end{theorem}
\begin{proof} {\it Case} (i).
A direct computation utilizing formulas \eqref{4.15}, 
\eqref{4.35} shows that
\begin{align}
 g(z,k_0)&= [M_-(z,k_0) -M_+(z,k_0)]^{-1}, 
\lb{4.36} \\ 
G(z,k_0,k_0+1)&=
-[M_-(z,k_0) -M_+(z,k_0)]^{-1}M_+(z,k_0)A(k_0)^{-1}. \lb{4.37}
\end{align}
In addition, the matrix $A(k_0)$ can be read off from the 
asymptotic expansion of $G(z,k_0,k_0+1,z)$ as $z\to\infty$ by a
simple Neumann series-type argument. Indeed, by \eqref{4.10} 
\begin{equation}
m_{\pm}(z,k_0)=-z^{-1}I_m -z^{-2}B(k_0)
-z^{-3}(A(k_0)^2+B(k_0)^2) +O(|z|^{-4})
\end{equation}
and hence by \eqref{4.11} and \eqref{4.12}
 \begin{align}
M_{+}(z,k_0)&
\underset{|z|\to\infty}{=}
 -z^{-1}A(k_0)^2 +O(|z|^{-2}), \lb{4.38} \\
M_{-}(z,k_0)&
\underset{|z|\to\infty}{=}
 -zI_m +B(k_0)+z^{-1}A(k_0)^2+O(|z|^{-2}) \lb{4.39}
\end{align}
along any ray with $\arg(z)=\pi -\varepsilon$, $0<\varepsilon<\pi$. 
Therefore, combining 
\eqref{4.36}, \eqref{4.37},  \eqref{4.38}, and \eqref{4.39}
one infers that
\begin{equation}\lb{4.40}
G(z,k_0,k_0+1)\underset{|z|\to\infty}{=}-z^{-2}A(k_0)
+O(|z|^{-3}), \quad k_0\in \bbZ
\end{equation}
along any ray with $\arg(z)=\pi-\varepsilon$, $0<\varepsilon<\pi$,
in complete analogy with the corresponding scalar case (see, e.g.,
\cite{Te98}, \cite[Sect.~6.1]{Te00}). Given $A(k_0)$ and using 
\eqref{4.36} one can easily solve \eqref{4.37} for $M_\pm (z,k_0)$,
\begin{align}
M_+(z,k_0)&=-g(z,k_0)^{-1}G(z,k_0,k_0+1)A(k_0), \lb{4.41} \\
M_-(z,k_0)&=M_+(z,k_0)+g(z,k_0)^{-1}\no \\
&=g(z,k_0)^{-1}(I_m-G(z,k_0,k_0+1)A(k_0)), \lb{4.42}
\end{align}
which in turn determines
$A(k)$, $B(k)$ for all $k\in\bbZ$ by Theorem~\ref{t4.2}. 

\noindent {\it Case} (ii). The data (ii) 
are perhaps the closest discrete analog of the 
continuous version in Theorem~\ref{t3.5}. As before (cf.
\eqref{4.40}), one can determine $A(k_0)$ from the asymptotic
expansion of $G(z,k_0,k_0+1)+G(z,k_0+1,k_0)$,
\begin{equation}\lb{4.44}
G(z,k_0,k_0+1)+G(z,k_0+1,k_0)\underset{|z|\to\infty}{=}
-2z^{-2}A(k_0)+O(|z|^{-3}), \quad k_0\in \bbZ.
\end{equation}
Next, using \eqref{4.15}, one observes that $M_+(z,k_0)$
is the unique solution of the following Sylvester-type matrix
equation
\begin{align}
&G(z,k_0+1,k_0) +G(z,k_0,k_0+1)
\no \\
&=-A(k_0)^{-1}M_+ (z,k_0)g(z,k_0)-g(z,k_0)M_+(z,k_0)A(k_0)^{-1}
 \lb{4.43} 
\end{align}
(cf.~\eqref{3.44}) by Lemmas~\ref{l3.3} and \ref{l3.4}.
Determining $M_-(z,k_0)$ from \eqref{4.42} completes the proof of
case (ii) by Theorem~\ref{t4.2} as in case (i).

\noindent {\it Case} (iii). Since $A(k_0)$ is  known, 
the data (iii), \eqref{4.15}, and \eqref{4.35} imply
\begin{equation}\lb{4.48} 
g(z,k_0+1) = A(k_0)^{-1}M_+ (z,k_0)[M_-(z,k_0) 
-M_+(z,k_0)]^{-1}M_-(z,k_0)A(k_0)^{-1}.
\end{equation} 
Together with formula \eqref{4.36} for $g(z,k_0)$ this leads to 
the following equation for $M_+(z,k_0)$,
\begin{equation}
M_+(z,k_0)+M_+(z,k_0)g(z,k_0)M_+(z,k_0)
=A(k_0)g(z,k_0+1)A(k_0). \lb{4.49}
\end{equation}
Using the asymptotics  \eqref{4.38}, \eqref{4.39}
and the representation
\begin{equation}\lb{4.50}
g(z,k)=-z^{-1}I_m+O(|z|^{-2}),\quad k\in \bbZ,
\end{equation}
 one infers
\begin{align} 
&\|M_+(z,k_0)\|_{\bbC^{m\times m}}<1,\quad
\|A(k_0)g(z,k_0+1)A(k_0)\|_{\bbC^{m\times m}}<\frac{1}{2}, \no \\
&\text{ and } \, \|g(z,k_0)\|_{\bbC^{m\times m}}<\frac{1}{2} 
\lb{4.51}
\end{align}
for $\Im(z)>0$ sufficiently large.
Thus, by Lemma~\ref{l.ric}, $M_+(z,k_0)$ is uniquely 
determined by the data (iii) for $\Im(z)>0$ large enough 
and hence for all $z\in \bbC_+$ by analytic continuation. The 
Weyl matrix $M_-(z,k_0)$ is then determined  by (cf. \eqref{4.42})
\begin{equation}\lb{4.52}
M_-(z,k_0)=M_+(z,k_0)+g(z,k_0)^{-1},
\end{equation}
completing the proof.
 \end{proof}

\begin{remark}
We note that $M_+ (z,k_0)$ can explicitly be represented by the
contour integral (cf. \eqref{3.44a})
\begin{align}
&M_+(z,k_0) \no \\
&=(2\pi i)^{-1}\oint_\Gamma d\zeta\,
(g(z,k_0)-\zeta A(k_0)^{-1})^{-1}
[G(z,k_0,k_0+1)+G(z,k_0+1,k_0)]\times \no \\
&\hspace*{2.1cm} \times (g(z,k_0)+\zeta A(k_0)^{-1})^{-1}.
\lb{4.46}
\end{align}
Here $\Gamma$ is a closed  clockwise oriented Jordan contour such 
that $\spec (A(k_0)g(z,k_0))$ has winding number $+1$ and 
$\spec(-g(z,k_0)A(k_0))$  has winding number $0$ with respect to
$\Gamma$. Without loss of generality we may assume that
$\Gamma\subset \bbC_+$. The existence of such a contour $\Gamma$
follows from the  observation that
\begin{equation}\lb{4.47}
\spec (g(z,k_0)A(k_0))=
\spec (A(k_0)^{1/2}g(z,k_0)A(k_0)^{1/2})
\subset \bbC_+, \quad z\in \bbC_+
\end{equation}
by Lemma~\ref{l3.4}.
\end{remark}

The special scalar case $m=1$ of Theorem~\ref{t4.4} is known 
and has recently been derived in \cite{Te98} (see also 
\cite[Sect.~2.7]{Te00}. Condition (iii)  in Theorem~\ref{t4.4} is
specific to the Jacobi case. In the corresponding Schr\"odinger 
case, the corresponding set of data does not even uniquely
determine the potential in the scalar case $m=1$ (see, 
e.g., \cite{GS96a}, \cite{Te98}).

The local analog of Theorem~\ref{t3.5a} in the discrete 
context then reads as follows.

\begin{theorem} \lb{t4.6}
Suppose $0< \varepsilon<\pi$. Let $H_j$, $j=1,2$ be two Jacobi
operators in 
$\ell^2 (\bbZ)^m$, with coefficients $A_j=\{A_{j,k}\}_{k\in\bbZ}$,
$B_j=\{B_{j,k}\}_{k\in\bbZ}$ satisfying Hypothesis~\ref{h4.1}, and
denote by $G_j(z,k,\ell)$, $g_j(z,k)=G_j(z,k,k)$ $(k,\ell\in \bbZ)$
the Green's matrices  associated with $H_j$, $j=1,2$. Moreover,
assume there exists an $N\in\bbN$, $N\geq 2$, such that  one of the
following conditions holds.

\noindent $(i)$
\begin{align}\lb{4.53} 
&\|g_1(z,k_0)-g_2(z,k_0)\|_{\bbC^{m\times m}} +
\|G_1(z,k_0,k_0+1)-G_2(z,k_0,k_0+1)\|_{\bbC^{m\times m}} \no \\
&\, \underset{|z|\to\infty}{=} O(|z|^{-(N+2)});
\end{align}

\noindent $(ii)$ 
\begin{align}\lb{4.54} 
&\|g_1(z,k_0)-g_2(z,k_0)\|_{\bbC^{m\times m}} +
\|G_1(z,k_0,k_0+1)-G_2(z,k_0,k_0+1) \no \\ &
+G_1(z,k_0+1,k_0)-G_2(z,k_0+1,k_0)\|_{\bbC^{m\times m}}
\underset{|z|\to\infty}{=} O(|z|^{-(N+2)});
\end{align}

\noindent $(iii)$ $A_1(k_0)=A_2(k_0)$ is given and 
\begin{equation}\lb{4.55} 
\|g_1(z,k_0)-g_2(z,k_0)\|_{\bbC^{m\times m}} +
\|g_1(z,k_0+1)-g_2(z,k_0+1)\|_{\bbC^{m\times m}}
\underset{|z|\to\infty}{=} O(|z|^{-(N)})
\end{equation}
along a ray with $\arg(z)=\pi -\varepsilon$. Then 
\begin{align}
A_{1}(k_0\pm k)&=A_{2}(k_0\pm k), 
\quad k=0, \, ...\, , \lfloor N/2\rfloor -1,
\no \\
 \quad B_{1}(k_0\pm k)&=B_{2}(k_0\pm k), \quad k=0, \, ...\, ,
 \lfloor (N-1)/2 \rfloor.
\lb{4.56}
\end{align}
\end{theorem}
\begin{proof} {\it Case} (i).   Since $N\ge 2$ in \eqref{4.53},
\eqref{4.40} implies that $A_1(k_0)=A_2(k_0)$ and we denote the
latter by $A(k_0)$. By \eqref{4.41} the $m\times m$
 Weyl-Titchmarsh matrices $M_{j,+}(z,k_0)$ 
associated with $H_j$, $j=1,2$ can be represented as
\begin{equation}\lb{4.57}
M_{j,+}(z,k_0)=-g_{j}(z,k_0)^{-1}G_{j}(z,k_0,k_0+1)A(k_0), 
\quad j=1,2 
\end{equation}
and they admit the estimate
\begin{align}
&\|M_{1,+}(z,k_0)-M_{2,+}(z,k_0)\|_{\bbC^{m\times m}} \no \\
&\le \|A(k_0)\|_{\bbC^{m\times m}} 
\|g_{1}(z,k_0)^{-1}G_{1}(z,k_0,k_0+1)-g_{2}(z,k_0)^{-1}G_{2}
(z,k_0,k_0+1)\|_{\bbC^{m\times m}} \no \\
&\le \|A(k_0)\|_{\bbC^{m\times m}} 
\|g_{1}(z,k_0)^{-1}\|_{\bbC^{m\times m}} 
\|G_{1}(z,k_0,k_0+1)-G_{2}(z,k_0,k_0+1)\|_{\bbC^{m\times m}} \no \\
&\quad + \|A(k_0)\|_{\bbC^{m\times m}} 
 \|g_{1}(z,k_0)^{-1}-g_{2}(z,k_0)^{-1}\|_{\bbC^{m\times m}}\,
\|G_{2}(z,k_0,k_0+1)\|_{\bbC^{m\times m}} \no \\
&\le\|A(k_0)\|_{\bbC^{m\times m}}
\|g_{1}(z,k_0)^{-1}\|_{\bbC^{m\times m}}  
\|G_{1}(z,k_0,k_0+1)-G_{2}(z,k_0,k_0+1)\|_{\bbC^{m\times m}}
\no \\ 
&\quad + \|A(k_0)\|_{\bbC^{m\times m}} 
\|g_{1}(z,k_0)^{-1}\|_{\bbC^{m\times m}}
\|g_{1}(z,k_0)-g_{2}(z,k_0)\|_{\bbC^{m\times m}}\times \no \\
& \quad \times \|g_{2}(z,k_0)^{-1}\|_{\bbC^{m\times m}}
\|G_{2}(z,k_0,k_0+1)\|_{\bbC^{m\times m}}. \lb{4.58}
\end{align}
Using the asymptotics 
\begin{align}
g_j(z,k_0)&\underset{|z|\to\infty}{=}-z^{-1}I_m
+O(|z|^{-2}), \quad j=1,2, \lb{4.59} \\
G_j(z,k_0,k_0+1)&\underset{|z|\to\infty}{=}-z^{-2}A(k_0)
+O(|z|^{-3}), \quad j=1,2, \lb{4.60}
\end{align}
the estimate \eqref{4.58} combined with \eqref{4.53} proves
\begin{equation}\lb{4.61}
\|M_{1,+}(z,k_0)-M_{2,+}(z,k_0)\|_{\bbC^{m\times m}}
\underset{|z|\to\infty}{=} O(|z|^{-(N+1)}).
\end{equation}
Since
\begin{equation}\lb{4.62}
M_{j,-}(z,k_0)=M_{j,+}(z,k_0)+g_j(z,k_0)^{-1}, \quad j=1,2,
\end{equation}
\eqref{4.53} and \eqref{4.61} yield
\begin{equation}\lb{4.63}
\|M_{1,-}(z,k_0)-M_{2,-}(z,k_0)\|_{\bbC^{m\times m}} 
\underset{|z|\to\infty}{=}
O(|z|^{-N}).
\end{equation}
\eqref{4.27} and \eqref{4.63} yield
\begin{equation}
B_1(k_0)=B_2(k_0). \lb{4.63a}
\end{equation} 
The asymptotic expansions \eqref{4.61} and \eqref{4.63} together 
with \eqref{4.63a} then prove \eqref{4.56} by Theorem~\ref{t4.2a}.

\noindent {\it Case} (ii). Since $N\geq 2$ in \eqref{4.54},
\eqref{4.40} again implies that $A_1(k_0)=A_2(k_0)$ and we denote
the latter by $A(k_0)$. Given $\varepsilon \in (0, \pi)$, let
$\gamma$ be a closed  clockwise oriented Jordan contour 
\begin{equation}\lb{4.64}
\gamma\subset \{ z\in \bbC \, |\, \Re(z)>0 \}
\end{equation}
 such 
that $\spec (A(k_0))$ 
 has winding number $+1$ and
\begin{equation}\lb{4.65}
e^{i\varepsilon}\gamma\subset \bbC_+.
\end{equation}
Next, given $z\in \bbC_+$, 
$ \arg(z)=\pi -\varepsilon$, introduce the contour $\Gamma(z)$ in
 the upper  half-plane by
\begin{equation}\lb{4.66}
\Gamma(z)=-z^{-1}\gamma.
\end{equation}
 Using \eqref{4.59} one concludes that
$\spec (A(k_0)^{1/2}g_j(z,k_0)A(k_0)^{1/2})$   
has winding number $+1$ with respect to the contour 
$\Gamma(z)$ if $z$ lies on the ray $\arg(z)=\pi -\varepsilon$ 
for $|z|$ sufficiently large. For such $z$ the matrix-valued
functions $M_{j,+}(z,k_0)$ admit the representation
(cf. \eqref{4.46})
\begin{align}
M_{j,+}(z,k_0)=(2\pi i)^{-1}\oint_{\Gamma(z)} & d\zeta\,
\widehat G_j(z,k_0, \zeta )
[G_j(z,k_0,k_0+1)+G_j(z,k_0+1,k_0)]\times \no \\  
&\times \widehat G_j(z,k_0,-\zeta) , \quad j=1,2, \lb{4.67}
\end{align}
where
\begin{align}
&\widehat G_j(z,k_0, \zeta)=A(k_0)^{1/2}(A(k_0)^{1/2}
g_j(z,k_0)A(k_0)^{1/2}-\zeta)^{-1}A(k_0)^{1/2}, \lb{4.68} \\
&\hspace*{7.4cm}  j=1,2, \,\, \zeta\in \Gamma(z).
\no 
\end{align}
Introducing
\begin{equation}\lb{4.69}
f_j(z,k_0)=g_j(z,k_0)+z^{-1}I_m, \quad j=1,2,
\end{equation}
one infers by \eqref{4.59} that 
\begin{equation}\lb{4.70}
f_j(z,k_0)\underset{|z|\to\infty}{=}O(|z|^{-2}), \quad j=1,2.
\end{equation}
The asymptotics
\begin{align}
& \sup_{\zeta\in \Gamma(z)}\|(-z^{-1}A(k_0)+A(k_0)^{1/2}
f_j(z,k_0) A(k_0)^{1/2}\mp\zeta)^{-1}\|_{\bbC^{m\times m}} \no \\
&=|z| \, \sup_{\zeta\in \Gamma(z)}\|(A(k_0)-zA(k_0)^{1/2}
f_j(z,k_0) A(k_0)^{1/2}\pm z\zeta)^{-1}\|_{\bbC^{m\times m}} \no \\
&=|z|\sup_{\zeta\in \gamma}\|(A(k_0)\mp 
\zeta)^{-1}\|_{\bbC^{m\times m}}\times \no \\
& \quad \,\,\sup_{\zeta\in \gamma}\|(I_m-z(A(k_0)\mp
\zeta)^{-1}A(k_0)^{1/2} f_j(z,k_0)
A(k_0)^{1/2})^{-1}\|_{\bbC^{m\times m}} \no \\
&\underset{|z|\to\infty}{=}O(|z|)
\lb{4.71}
\end{align}
along the ray $\arg(z)=\pi -\varepsilon$, and the representation
\begin{align}
&\widehat G_j(z,k_0, \zeta)=A(k_0)^{1/2}(-z^{-1}A(k_0)+A(k_0)^{1/2}
f_j(z,k_0)
A(k_0)^{1/2}-\zeta)^{-1}A(k_0)^{1/2}, \no \\
&\hspace*{10cm}  j=1,2 \lb{4.72}
\end{align} 
then yield
\begin{equation}\lb{4.73}
\|\widehat G_j(z,k_0, \pm\zeta)
\|_{\bbC^{m\times m}}\underset{|z|\to\infty}{=}
O(|z|), \quad \zeta\in \Gamma(z),\quad j=1,2,
\end{equation}
along the ray $\arg(z)=\pi -\varepsilon$, uniformly with respect
to $\zeta \in \Gamma(z)$. As a consequence of \eqref{4.73} and
hypothesis \eqref{4.54} we have 
\begin{equation}\lb{4.74}
\|\widehat G_1(z,k_0, \pm\zeta)-\widehat 
G_2(z,k_0, \pm\zeta)\|_{\bbC^{m\times m}}
\underset{|z|\to\infty}{=}
O(|z|^{-N}), \quad \zeta\in \Gamma(z) 
\end{equation}
along the ray $\arg(z)=\pi -\varepsilon$, uniformly in $\zeta 
\in \Gamma(z)$.
Using \eqref{4.67}, \eqref{4.73}, and \eqref{4.74}
one concludes that
\begin{equation}\lb{4.75}
\|M_{1,+}(z,k_0)-M_{2,+}(z,k_0)\|_{\bbC^{m\times m}}
\underset{|z|\to\infty}{=}
O(|z|^{-(N+1)})
\end{equation}
and hence by \eqref{4.56} and \eqref{4.54} that 
\begin{equation}\lb{4.76}
\|M_{1,-}(z,k_0)-M_{2,-}(z,k_0)\|_{\bbC^{m\times m}}
\underset{|z|\to\infty}{=}
O(|z|^{-N})
\end{equation}
along the ray $\arg(z)=\pi -\varepsilon$. This proves \eqref{4.56}
by Theorem~\ref{t4.2a}.\\
\noindent {\it Case} (iii). One observes that 
 $M_{j,+}(z,k_0)$ satisfy (c.f. \eqref{4.49})
\begin{equation}
M_{j,+}(z,k_0)+M_{j,+}(z,k_0)g_j(z,k_0)M_{j,+}(z,k_0)
=A(k_0)g_j(z,k_0+1)A(k_0), \,\,\, j=1,2 \lb{4.77}.
\end{equation}
Taking into account the asymptotic expansions
\eqref{4.59} and \eqref{4.60},  hypothesis \eqref{4.55} implies
\eqref{4.75} and \eqref{4.76} by Corollary~\ref{cor}. This proves
\eqref{4.56} applying Theorem~\ref{t4.2a} again.
\end{proof}

Without going into further details, we note that the methods 
employed in this section can easily be adapted to the case of
uniformly bounded operator coefficients for Jacobi operators 
(i.e., for some $C>0$ and all $k\in\bbZ$, 
$\|A(k)\|_{\calB(\calH)}+\|B(k)\|_{\calB(\calH)}\leq C$,
$A(k)^{-1}\in\calB(\calH)$, $0<A(k)$, $B(k)=B(k)^*$,
 for some complex separable Hilbert space $\calH$).   

\vspace*{2mm}
\noindent {\bf Acknowledgements.} We are indebted to Steve 
Clark, Don Hinton, Boris Levitan, Mark Malamud, 
Fedor Rofe-Beketov, Alexander Sakhnovich, Lev Sakhnovich, and Barry
Simon for discussions and valuable hints regarding the literature.
F.G. and A.K. thank C.~Peck and T.~Tombrello  for the
hospitality of Caltech during the  final stages of this work.


\end{document}